\documentclass[a4paper,12pt]{article}
\usepackage{amssymb,amsmath}
\newcommand{\dz}{\delta_0}
\newcommand{\dw}{\delta_1}

\newcommand{\cl}[1]{\mathcal{#1}}
\newcommand{\card}{\#}

\newcommand{\cA}{\cl{A}}
\newcommand{\cB}{\cl{B}}

\newcommand{\bs}{\mathbf{S}}
\newcommand{\bso}{\mathbf{S}^{o(1)}}
\newcommand{\Z}{\mathbb{Z}}

\newcommand{\twosum}[2]{\sum_{\substack{#1\\#2}}}
\newcommand{\ep}{\varepsilon}
\newcommand{\beql}[1]{\begin{equation}\label{#1}}
\newcommand{\eeq}{\end{equation}}

\newcommand{\andd}{\;\;\;\mbox{and}\;\;\;}

\newtheorem{theorem}{Theorem}
\newtheorem{Prop}{Proposition}

\newtheorem{lemma}{Lemma}

\begin{document}

\title{The Differences Between Consecutive\\ Primes. V}

\author{D.R. Heath-Brown\\Mathematical Institute, Oxford}

\date{}

\maketitle


\section{Introduction}

In this paper we shall continue our investigations into the sum
\beql{bs}
\twosum{p_n\le x}{p_{n+1}-p_n\ge\sqrt{p_n}}(p_{n+1}-p_n).
\eeq
It was shown by Wolke \cite{W} that the sum is $O(x^{1-\delta})$ (with
$\delta=1/30$), thereby
answering a question
of Erd\H{o}s. The exponent was improved firstly to $5/6+\ep$ for
any fixed $\ep>0$, and then to $3/4+\ep$, by the author \cite{HBI},
\cite{HBIII}. A further reduction, to $25/36+\ep$, was achieved by
Peck \cite{Peck}, and the present record is held by Matom\"{a}ki
\cite{Mat}, with exponent $2/3$.

One would conjecture that the sum (\ref{bs}) only contains at most
the terms with $p_n=3$, 7, 13, 23, 31 and 113, and hence is bounded.
However we are far from proving
this, even under the Riemann Hypothesis. The latter assumption allows
for an estimate $O(x^{1/2}(\log x)^2)$, as was proved by Selberg
\cite{Selb}.  The Lindel\"{o}f Hypothesis similarly implies that the
sum is $O_{\ep}(x^{1/2+\ep})$, for any positive $\ep$, as was shown
in the fourth paper of this series, \cite{RHBIV}.

Our goal is to improve the unconditional estimates as follows.
\begin{theorem}\label{T1}
For any fixed $\ep>0$ we have
\[\twosum{p_n\le x}{p_{n+1}-p_n\ge\sqrt{p_n}}(p_{n+1}-p_n)\ll_{\ep} x^{3/5+\ep}.\]
\end{theorem}
We should view the exponents $2/3$ and $3/5$ as being
\[1/2+1/6\;\;\;\mbox{and}\;\;\; 1/2+1/10\]
respectively, so that we have reduced the excess over $1/2$
by $40\%$, from $1/6$ down to $1/10$. (For comparison, Heath-Brown
\cite{HBI} gives roughly a $29\%$ improvement over Wolke \cite{W};
Heath-Brown \cite{HBIII} sharpens \cite{HBI} by $25\%$; Peck \cite{Peck}
improves on \cite{HBIII} by about $22\%$; and Matom\"{a}ki reduces the
excess in Peck's exponent by some $14\%$.)

In fact we prove a stronger result than Theorem \ref{T1}.
\begin{theorem}\label{T2}
  For any fixed $\ep>0$ the measure of the set of $y\in[0,x]$ such that
  \[    \max_{0\le h\le\sqrt{y}}
  \left|\pi(y+h)-\pi(y)-\int_y^{y+h}\frac{dt}{\log t}\right|\ge
  \frac{\sqrt{y}}{(\log y)(\log\log y)} \]
is $O_{\ep} (x^{3/5+\ep})$.
\end{theorem}
If $p_n\le y\le p_{n+1}-\tfrac12\sqrt{p_{n+1}}$ then
$\pi(y+\tfrac12\sqrt{y})-\pi(y)=0$. Thus each $p_n\le x$ for which
$p_{n+1}-p_n\ge\sqrt{p_n}$ 
contributes an interval whose length is at least
$\tfrac13(p_{n+1}-p_n)$, say, to the set in Theorem \ref{T2}, provided
that $p_n$ is large enough.  Thus Theorem~\ref{T1} is a corollary of
Theorem \ref{T2}.

We remark that the analysis in this paper would be very considerably
simplified if our question had been about gaps $p_{n+1}-p_n$ of size
at least $p_n^{1/2+\ep}$, rather than $p_n^{1/2}$.
There would be no ``bad'' ranges to be handled
by sieve upper bounds, and one could merely have used the
generalized Vaughan identity, rather than the Buchstab
formula. The situation is analogous to that in the author's papers
\cite{hbi} and \cite{hbg}, proving $\pi(x+y)-\pi(x)\sim y/\log x$ for
Huxley's range $x^{7/12+\ep}\le y\le x$, and for $x^{7/12}\le y\le x$
respectively. 

With the exception of Matom\"{a}ki's work, all previous
results on the sum (\ref{bs}) could have been adapted
to prove a corresponding
version of Theorem~\ref{T2}. Matom\"{a}ki uses sieve methods in an
essential way, so that her method shows the sparsity of values $y$
where $\pi(y+\sqrt{y})-\pi(y)\le cy/\log y$ for some small positive
constant $c$.  In contrast, our approach only uses sieve methods to
handle relatively minor contributions to
$\pi(y+\sqrt{y})-\pi(y)$. Theorem~\ref{T1} could undoubtedly
be improved further by deploying sieve methods in the same way that
Matom\"{a}ki does. We have decided against doing this largely from
laziness, but partly so as to demonstrate more clearly the power of
our primary new tool, Proposition \ref{P1}, described below.

Our approach to the sum (\ref{bs}) uses the standard mean and large values
estimates for Dirichlet polynomials, which arise naturally in this
context when one applies a sieve decomposition to the problem. The
Dirichlet polynomials one encounters are typically products of shorter
polynomials of the form $\sum_{N<n\le 2N}p^{-s}$, the sum being over
primes.

The mean value theorem for Dirichlet polynomials (Montgomery
\cite[Theorem 6.1]{mont}) shows that
\beql{mvdp}\int_0^T\left|\sum_{m\le M}a_m m^{-it}\right|^2 dt\ll
(T+M)\sum_{m\le M}|a_m|^2.
\eeq
This is quite efficient when $M\ll T$ and the coefficients $a_m$ are
fairly even in size.  Our main new tool is the following quite different mean
value estimate, which remains useful for certain longer Dirichlet
polynomials.
\begin{Prop}\label{P1}
Let $T\ge 1$ and let $\cl{M}$ be a set of
distinct integers $m$ in $(0,T]$, of cardinality $\#\cl{M}\le R$, with
  associated complex coefficients $\zeta_m$ of modulus 1.
Suppose we are given a positive integer $N$ and complex coefficients
$q_1,\ldots,q_N$.  Then we have
\begin{eqnarray*}
\lefteqn{\int_0^T\left|\sum_{m\in\cl{M}}\zeta_m m^{-it}\right|^2
\left|\sum_{n\le N}q_n n^{-it}\right|^2dt}\\
&\ll_{\ep}&\left(N^2R^2+(NT)^{\ep}\{NRT+NR^{7/4}T^{3/4}\}\right)\max_n|q_n|^2
\end{eqnarray*}
for any fixed $\ep>0$.  
  \end{Prop}

For the proof we refer the reader to the author's paper
\cite{smooth}. To see the strength of this result we observe that the
term $N^2R^2$ corresponds to the (square of the) maximum value that
the product of our two Dirichlet polynomials could attain, while the term
$NRT$ is what one would get if one had square root cancellation
throughout the range $[0,T]$.  Thus the bound is sub-optimal largely
because of the term $NR^{7/4}T^{3/4}$. When $R\le T^{1/3}$ one has
$NR^{7/4}T^{3/4}\le NRT$, so that our result is essentially best
possible in this case.

Estimates of the type in Proposition \ref{P1} originate with the work
of Yu \cite{YG}, who gave a bound $O_{\ep}((N^2R^2+NRT)(NT)^{\ep})$
subject to the Lindel\"{o}f hypothesis, and used it to show that
\beql{yg1}
\sum_{p_n\le x}(p_{n+1}-p_n)^2\ll_{\ep}x^{1+\ep}
\eeq
under the same assumption. Since we only obtain an optimal estimate in
Proposition \ref{P1} when $R\le T^{1/3}$ it turns out that we are
unable to say anything useful about gaps $p_{n+1}-p_n$ shorter than
$p_n^{1/3}$. This does not preclude a new unconditional result for the sum in
(\ref{yg1}), but in the present paper we will restrict our attention
to gaps of size $p_n^{1/2}$ or more.

{\bf Acknowledgements.}  This work was partly supported by EPSRC grant number
EP/K021132X/1.  For a further part of the preparation of this  paper the
author was
supported by the NSF under Grant No.\ DMS-1440140,  while  in residence at 
the {\em  Mathematical Sciences Research Institute} in Berkeley, California,
during the Spring 2017 semester.

\section{Structure of the Proof, and the Choice of Parameters}

Theorem \ref{T2} compares the number of primes in an interval
$(y,y+h]$ with its expected main term, for variable $h$.  Our argument
  begins with some preliminary steps to replace the interval $(y,y+h]$ by
$(y,y+y\dz]$ for a suitably small $\dz>0$ independent of $y$. Rather than
compare the number of primes in $(y,y+y\dz]$ with its expected main
  term, we find it convenient to compare with the number of primes in
  a longer interval $(y,y+y\dw]$.

The argument then goes on to apply sieve methods to both $(y,y+y\dz]$
  and $(y,y+y\dw]$.  This has two effects. Firstly it enables us to
    remove certain short ranges of variables that would otherwise be
    awkward to handle. Secondly it allows us to translate the problem
    into one involving products of Dirichlet polynomials. At this
    point we use the key idea from Yu \cite{YG}, coding the points $y$
    for which $(y,y+y\dz]$ does not have the expected number of primes
      into a Dirichlet polynomial.  This process produces Proposition
      \ref{P2} in Section \ref{secp2},
      which is a major waypoint in our argument.

Proposition \ref{P2} requires us to estimate the mean value of a
product of Dirichlet polynomials, and we do this using a variety of
well established techniques, in combination with Proposition
\ref{P1}. In particular we use Vinogradov's zero-free region, various
forms of the ``Large Values'' estimate for Dirichlet polynomials, and
the classical mean value estimate (\ref{mvdp}). This stage of the
argument requires us to examine several separate cases. For these estimates to
produce a suitable saving it is crucial that certain critical ranges
for the lengths of the Dirichlet polynomials are avoided, and these
are the ranges that the initial sieve argument eliminates.

The ranges we will avoid take the
shape $[x^{1/\ell-\eta},x^{1/\ell+\eta}]$ for certain positive
integers $\ell$, and by removing these we are able to make savings of
factors of order $x^{c\eta}$ for certain constants $c>0$. Here
$\eta=\eta(x)$ is a small
function of $x$ which we will specify in a moment.

At other points in the argument the saving we obtain is related to the
available zero-free region for the Riemann zeta-function, which allows
us to improve on the trivial bound by factors of the type
$\exp((\log x)^{\theta})$ for certain constants $\theta\in(0,1)$.
With this in mind we define
\beql{bsd}
\bs=\bs(x)=\exp\left((\log\log x)^{11}\right)
\eeq
with a view to saving at least a positive power of $\bs$ in the
various key arguments.  This means, conversely, that we can afford to
loose factors $\bso$, since they will be more than compensated for by
the gain of a power of $\bs$.

We set
\beql{nubnd}
\nu=\nu(x)=(\log\log x)^5,
\eeq
\beql{z1nu}
z_1=z_1(x)=(4x)^{1/\nu}=x^{o(1)},
\eeq
\beql{etanu}
\eta=\eta(x)=(\log\log x)^{-12000},
\eeq
\beql{vpd}
\varpi=\varpi(x)=(\log\log x)^2,
\eeq
and
\beql{vpd1}
z_2=z_2(x)=z_1^{\varpi}=x^{o(1)}.
\eeq
Then
\beql{80}
  \log x\ll\bso,\;\;\;\;\;\;\mbox{and}\;\;\;\nu^{\nu}\ll\bso,
  \eeq
while
\beql{81}
\bs\ll\exp\left((\log x)^{\theta}\right)\ll x^{o(\eta)}\ll z_1 \;\;\;
\mbox{for any}\; \theta\in(0,1).
\eeq
We shall use all these bounds repeatedly in our argument.

For the entirety of the paper we will assume without further
comment that $x$ is sufficiently large.

\section{Preliminary Steps}

To prove Theorem \ref{T2} it suffices to establish the
corresponding bound when $y$ varies over a dyadic range
$(x,2x]$. 
  Theorem \ref{T2} requires us to estimate the measure of a
set of real numbers $y$, defined using the maximum over intervals
$(y,y+h]$ for varying $h$.  We begin by showing how to replace the
real variable $y$ by an integer variable $m$, and how to use an interval
whose length is a fixed fraction $\dz$ of its left-hand endpoint.
It will be convenient to write $\cl{I}(x)$ for
the set of $y\in(x,2x]$ such that
  \[\max_{0\le h\le\sqrt{y}}
  \left|\pi(y+h)-\pi(y)-\int_y^{y+h}\frac{dt}{\log t}\right|\ge
  \frac{\sqrt{y}}{(\log y)(\log\log y)},\]
 so that our goal is to estimate ${\rm Meas}(\cl{I}(x))$.
\begin{lemma}\label{L1}
  Let
\beql{dzch}
\dz=x^{-1/2}(\log\log x)^{-2}
\eeq
and
  \beql{Hch}
  H=x^{1/2}(\log\log x)^{-4}.
  \eeq
  Then there is a set of $R_0$
  distinct integers $m_1,\ldots, m_{R_0}\in[x/H,3x/H]$ for which
\[\left|\pi(mH(1+\dz))-\pi(mH)-
  \int_{mH}^{mH(1+\dz)}\frac{dt}{\log t}\right|
\ge \frac{\dz x}{3(\log x)(\log\log x)},\]
and such that
\[{\rm Meas}(\cl{I}(x))\ll  x^{1/2}R_0.\]
\end{lemma}

If $0\le h\le\sqrt{y}$  then the interval $(y,y+h]$
  is a union of at most
$\dz^{-1}y^{-1/2}$ disjoint subintervals $(y_1,y_1(1+\dz)]$, together
with a shorter interval at the end, of length $O(\dz x)$.
By the Brun--Titchmarsh theorem this last interval contains $O(\dz x/\log x)$
primes.  It follows that
\begin{eqnarray*}
  \left|\pi(y+h)-\pi(y)-\int_y^{y+h}\frac{dt}{\log t}\right|\hspace{6cm}\\
  \le \dz^{-1}y^{-1/2}\left|\pi(y_1(1+\dz))-\pi(y_1)-
  \int_{y_1}^{y_1(1+\dz)}\frac{dt}{\log t}\right| +O(\dz x/\log x)
\end{eqnarray*}
for some $y_1$ with $y\le y_1\le y+\sqrt{y}$.  Thus if $x<y\le 2x$ and
\beql{cnd}
\left|\pi(y+h)-\pi(y)-\int_y^{y+h}\frac{dt}{\log t}\right|\ge
\frac{\sqrt{y}}{(\log y)(\log\log y)}
\eeq
then there is some $y_1$ as above with
\begin{eqnarray*}
\left|\pi(y_1(1+\dz))-\pi(y_1)-
  \int_{y_1}^{y_1(1+\dz)}\frac{dt}{\log t}\right|
& \ge & \frac{\dz y}{2(\log y)(\log\log y)}\\
& \ge & \frac{\dz x}{2(\log x)(\log\log x)}.
\end{eqnarray*}
For each value $y$ satisfying (\ref{cnd}) we choose the smallest such
$y_1$ and define $m=1+[y_1/H]$, so that $x/H<m\le 3x/H$. In this way
we produce a collection of distinct integers $m_1,\ldots,m_{R_0} $.
Each such integer $m_i$ may arise from a range of
values for $y$, satisfying $y=m_iH+O(\sqrt{x})$. It therefore follows that
the measure we have to estimate in Lemma \ref{L1} is $O(x^{1/2}R_0)$, as
required. Moreover, if $m=1+[y_1/H]$ then $\pi(y_1(1+\dz))$  differs
from $\pi(mH(1+\dz))$  by
$O(H/\log x)$, by the Brun--Titchmarsh Theorem. Similarly
$\pi(y_1)=\pi(mH)+O(H/\log x)$. On the other hand,
\[\int_{y_1}^{y_1(1+\dz)}\frac{dt}{\log t}\]
differs from the corresponding integral between $mH$ 
and $mH(1+\dz)$
by $O(H/\log x)$.  Since $H/\log x=o(\dz x/(\log x)(\log\log x))$ we
therefore have
\begin{eqnarray*}
\left|\pi(mH(1+\dz))-\pi(mH)-
  \int_{mH}^{mH(1+\dz)}\frac{dt}{\log t}\right|
& \ge & \frac{\dz x}{3(\log x)(\log\log x)},
\end{eqnarray*}
for large enough $x$, and the lemma follows.
\medskip

For each $m\in[x/H,3x/H]$ we will locate the primes in the interval
\[\cA:=\cA(m)=\{n\in\Z: mH<n\le mH(1+\dz)\}\]
by sieving.  Rather than compare the number of primes with the
integral
\[\int_{mH}^{mH(1+\dz)}\frac{dt}{\log t}\]
it will be more convenient to work with the number of primes in a
long interval
\[\cB:=\cB(m)=\{n\in\Z: mH<n\le mH(1+\dw)\},\]
where
\beql{H0ch}
\dw:=\exp\{-(\log x)^{1/2}\}.
\eeq
We write $\cA^{(k)}=\cA^{(k)}(m)$ for the weighted set in which
$n\in\cA(m)$ has weight
$\tau_k(n)$, and similarly for $\cB^{(k)}$.  In what follows we will
re-number the integers $m_i$ in Lemma \ref{L1} as necessary.
We proceed to show the following.
\begin{lemma}\label{lf}
  There are integers $k=1,2$ or $3$ and
  $m_1,\ldots,m_{R_1}\in [x/H,3x/H]$ with $R_1\ge R_0/3$,
  satisfying
  \[\left|S(\cA^{(k)}(m_i),(4x)^{1/4})-
\dz^{-1}\dw S(\cB^{(k)}(m_i),(4x)^{1/4})\right|
\ge \frac{\dz x}{25(\log x)(\log\log x)}.\]
\end{lemma}

We begin the proof by using the prime number theorem with Vinogradov's
error term, whence
\begin{eqnarray*}
\pi(mH(1+\dw))-\pi(mH)&=&
\int_{mH}^{mH(1+\dw)}\frac{dt}{\log t}+O(\dw x(\log x)^{-2})\\
&=&
\dw\dz^{-1}\int_{mH}^{mH(1+\dz)}\frac{dt}{\log t}+O(\dw x(\log x)^{-2}),
\end{eqnarray*}
so that if $x$ is large enough we have
\[\left|\pi(mH(1+\dz))-\pi(mH)-\dz\dw^{-1}\{\pi(mH(1+\dw))-\pi(mH)\}\right|\]
\[\ge \frac{\dz x}{4(\log x)(\log\log x)}\]
for $m=m_1,\ldots,m_{R_0}$.

We now consider the weighted set $\cA_*$, consisting of 
elements $n\in\cA$ weighted by
\[3-\tfrac32 \tau(n)+\tfrac13 \tau_3(n).\]
One may check that this takes the
value 1 when $n$ is prime, and that it vanishes for square-free $n$
having 2 or 3 prime factors.  It follows that $S(\cA_*,(4x)^{1/4})$
differs from
\[\pi(mH(1+\dz))-\pi(mH)\]
by the contribution from integers which have a factor 
$p^2\in[(4x)^{1/2},4x]$. When
$(4x)^{1/4}\le p\le x^{1/3}$ there are $O(1+x\dz p^{-2})$ multiples
of $p^2$, producing 
total contribution $O(x^{1/3})$. On the other hand, if $p^2t\in\cA$ 
with $p\ge x^{1/3}$, then $t\le 4x^{1/3}$.  For each such $t$ the
prime $p$ is restricted to an interval $(\sqrt{mH/t},\sqrt{(mH(1+\dz))/t}]$
  of length $O(1)$, so that the total contribution in this case is also
  $O(x^{1/3})$.  It follows that
 \[S(\cA_*,(4x)^{1/4})=\pi(mH(1+\dz))-\pi(mH)
 +o(\frac{\dz x}{(\log x)(\log\log x)}).\]
 We define $\cB_*$ analogously, and find this time that the error is
 \[\ll\sum_{(4x)^{1/4}\le p\le(4x)^{1/2}}(1+x\dw p^{-2})\ll x^{3/4}
=o(\frac{\dw x}{(\log x)(\log\log x)}).\]
  We may then deduce that
\[\left|S(\cA_*,(4x)^{1/4})-\dz\dw^{-1}S(\cB_*,(4x)^{1/4})\right|
\ge \frac{\dz x}{5(\log x)(\log\log x)},\]
for large enough $x$. Lemma \ref{lf} then follows.

The integer $k$ appearing in Lemma \ref{lf}
will be fixed for the rest of the proof.

\section{The First Sieve Stage}\label{S1}

In this section we introduce our first sieve process, and show that
terms in which certain variables lie in awkward ranges make a
negligible contribution to $S(\cA^{(k)}(m_i),(4x)^{1/4})$ and
$S(\cB^{(k)}(m_i),(4x)^{1/4})$.

We begin by noting that
$S(\cA^{(k)}(m),(4x)^{1/4})$ counts products $n_1\ldots n_k$ in $\cA(m)$
for which each factor $n_i$ has no prime divisor below $(4x)^{1/4}$.
We will use the parameters $\nu$ and $z_1$ given by (\ref{nubnd}) and
(\ref{z1nu}). We now define
\[\Pi_1:=\prod_{p<z_1}p,\andd \Pi_2=\prod_{z_1\le p<(4x)^{1/4}}p,\]
so that
\[S(\cA^{(k)}(m),(4x)^{1/4})=\sum_{q_1,\ldots, q_k\mid \Pi_2}
\mu(q_1)\ldots\mu(q_k)N_k(\cA,q_1\ldots q_k),\]
where
\beql{Ndef}
N_k(\cA,q):=\card\{h_1,\ldots,h_k:
qh_1\ldots h_k\in\cA(m),\,(h_1\ldots h_k,\Pi_1)=1\},
\eeq
and similarly for $\cB$. Each $q_i$ is composed of various prime
factors $p$, which belong to dyadic intervals of the type $(2^s,2^{s+1}]$.
  Similarly, we can decompose the range for the variables $h_i$ into
  dyadic intervals.

  We proceed to show that there is a negligible contribution from
  terms with a divisor close to a reciprocal
  power $x^{1/\ell}$, say. 
  \begin{lemma}\label{NL}
    Let $\eta$ be given by (\ref{etanu})
    and let $\cA_{\dagger}(m)$ be the set of integers in $\cA(m)$
    having a divisor in the range $[x^{1/\ell-2\eta},x^{1/\ell+2\eta}]$
    for some integer $\ell\in [4,\nu+2]$.  Then the number of
    $2k$-tuples $(q_1,\ldots,q_k,h_1,\ldots,h_k)$ with
    $q_1\ldots q_kh_1\ldots h_k$ in $\cA_{\dagger}(m)$, and such that
$q_1,\ldots, q_k\mid \Pi_2$ and $(h_1\ldots h_k,\Pi_1)=1$ is
\[\ll\dz x(\log x)^{-1}(\log\log x)^{-3/2}.\]
If we define $\cB_{\dagger}$ similarly then we get an analogous bound
with $\dw$ in place of $\dz$.
  \end{lemma}
If the prime divisors $p_i$ of $q_1\ldots q_k$ lie in dyadic intervals
$I_1,\ldots,I_t$, say, and the $h_i$ lie in dyadic intervals
$J_1,\ldots,J_k$, then $t+k\le k\nu+k$.  The lemma then shows that
there is a negligible contribution from
those collections of intervals for which any product from
$I_1,\ldots,I_t,J_1,\ldots,J_k$ lies in
$[x^{1/\ell-\eta},x^{1/\ell+\eta}]$. Here we use the fact that
$2^{k\nu+k}\le x^{\eta}$ for large $x$, by (\ref{nubnd}) and
(\ref{etanu}). In particular the lemma shows that we can reduce the
sieving range, restricting our attention to divisors $q$ of $\Pi_2$
whose prime factors satisfy $p<x^{1/4-\eta}$.

For the proof of the lemma we begin by observing that an integer
$n\in\cA$ arises
in at most $\tau_6(n)\le\tau(n)^5$ ways as
$n=h_1q_1\ldots h_kq_k$, so that for each $\ell$
the contribution we have to consider
is
\[\ll \sum_{x^{1/\ell-2\eta}\le d\le x^{1/\ell+2\eta}}
\twosum{n\in\cA}{(n,\Pi_1)=1,\,d\mid n}\tau(n)^5.\]
By Cauchy's inequality this is at most $\Sigma_1^{1/2}\Sigma_2^{1/2}$,
with
\[\Sigma_1=\sum_d\twosum{n\in\cA}{(n,\Pi_1)=1,\,d\mid n}\tau(n)^{10}
\le\sum_{n\in\cA:\, (n,\Pi_1)=1}\tau(n)^{11}\]
  and
\[\Sigma_2=\sum_{x^{1/\ell-2\eta}\le d\le x^{1/\ell+2\eta}}
\twosum{n\in\cA}{(n,\Pi_1)=1,\,d\mid n}1.\]
We now apply the following lemma, which is an immediate corollary of the
theorem of Shiu \cite{shiu}.
\begin{lemma}\label{lsh}
Suppose that $X,Y,z\ge 2$ are real numbers such that 
$X^c\le Y\le X$, for some constant $c>0$, and let $N$ be a positive
integer.  Then
\[\twosum{X<n\le X+Y}{p\mid n\implies p\ge z}\tau(n)^N\ll_{c,N}
\frac{Y}{\log X}\left(\frac{\log X}{\log z}\right)^{2^N}.\]
\end{lemma}
This produces
  \[\Sigma_1  \ll \frac{\dz x}{\log x}\nu^{2^{11}}.\]
Moreover, since $d\le x^{1/\ell+2\eta}\le x^{1/3}$, a
simple sieve upper bound yields
\[\twosum{n\in\cA}{(n,\Pi_1)=1,\,d\mid n}1\ll \frac{\dz x}{d\log z_1},\]
  and
\[\twosum{x^{1/\ell-2\eta}\ll d\ll x^{1/\ell+2\eta}}{(d,\Pi_1)=1}d^{-1}\ll
\frac{\eta\log x}{\log z_1},\]
whence
\[\Sigma_2\ll \frac{\dz x\eta(\log x)}{(\log z_1)^2}.\]
We therefore conclude that each value for $\ell$ contributes
\[\ll (\Sigma_1 \Sigma_2)^{1/2}\ll \frac{\dz x}{\log z_1}\nu^{2^{10}}\eta^{1/2}.\]
Since $\ell\ll\nu$ the total is
\[\ll\dz x(\log x)^{-1}\nu^{1026}\eta^{1/2}\ll
\dz x(\log x)^{-1}(\log\log x)^{-3/2},\]
in view of (\ref{nubnd}) and
(\ref{etanu}).  This proves Lemma \ref{NL} for $\cl{A}$, and the
treatment of $\cB$ is similar.
\medskip

Other terms we wish to remove are those in which $q_1\ldots q_k$ has
two or more prime factors $p_1,p_2$ lying in the same dyadic interval
$(2^s,2^{s+1}]$. 
\begin{lemma}\label{NL2}
  We have
\[\twosum{z_1\le p_1,p_2<x^{1/4-\eta}}{p_1/2\le p_2\le p_1}\;
\twosum{q_1,\ldots, q_k\mid \Pi_2}{p_1p_2\mid q_1\ldots q_k}
N_k(\cA,q_1\ldots q_k)\ll \dz x(\log x)^{-1}(\log\log x)^{-3/2},\]
and
\[\twosum{z_1\le p_1,p_2<x^{1/4-\eta}}{p_1/2\le p_2\le p_1}\;
\twosum{q_1,\ldots, q_k\mid \Pi_2}{p_1p_2\mid q_1\ldots q_k}
N_k(\cB,q_1\ldots q_k)\ll \dw x(\log x)^{-1}(\log\log x)^{-3/2}.\]
\end{lemma}

Arguing as for Lemma \ref{NL} we see that the total contribution
from the first sum is
\[\twosum{z_1\le p_1,p_2<x^{1/4-\eta}}{p_1/2\le p_2\le p_1}
\twosum{n\in\cA}{(n,\Pi_1)=1,\,p_1p_2\mid n}\tau_6(n).\]
This time Cauchy's inequality gives a bound
$\ll(\Sigma_1\Sigma_3)^{1/2}$, with $\Sigma_1$ as before, and
\[\Sigma_3=\twosum{z_1\le p_1,p_2<x^{1/4-\eta}}{p_1/2\le p_2\le p_1}
\twosum{n\in\cA}{(n,\Pi_1)=1,\,p_1p_2\mid n}1.\]
A simple sieve upper bound shows that the
inner sum above is
\[\ll\frac{\dz x/p_1p_2}{\log \min(z_1,\dz x/p_1p_2)}
\ll\frac{\dz x/p_1p_2}{\log z_1}+
\frac{\dz x/p_1p_2}{\log\dz x/p_1p_2},\]
whence
\begin{eqnarray*}
\Sigma_3&\ll&
\sum_{p_1}\left\{\frac{\dz x}{p_1(\log p_1)(\log z_1)}+
\frac{\dz x}{p_1(\log p_1)(\log\dz x/p_1^2)}\right\}\\
&\ll& \frac{\dz x}{(\log z_1)^2}+\frac{\dz x\log\log x}{(\log x)^2}\\
&    \ll& \frac{\dz x}{(\log z_1)^2},
\end{eqnarray*}
by (\ref{nubnd}) and (\ref{z1nu}).
 It then follows that the total contribution from terms where two
 prime factors lie in the same dyadic interval is
 \[\ll (\Sigma_1\Sigma_3)^{1/2}\ll
 \left(\frac{\dz x}{\log x}\nu^{2^{11}}\right)^{1/2}
 \left(\frac{\dz x}{(\log z_1)^2}\right)^{1/2}
 \ll \frac{\dz x}{(\log x)^{3/2}}\nu^{2^{10}+1}.\]
By the choice of $\nu$ in (\ref{nubnd}) this will be
suitably small, which completes the treatment of $\cA$.  The proof for 
$\cB$ is similar.

\section{A Second Sieve Operation}

The variables $h$ in (\ref{Ndef}) are now constrained to lie in
dyadic ranges which avoid the intervals
$[x^{1/\ell-\eta},x^{1/\ell+\eta}]$. We now write $\xi(h)=1$ if 
$(h,\Pi_1)=1$, and $\xi(h)=0$ otherwise.  This
is satisfactory when $h\le x^{1/4}$, but for larger $h$ we
need to pick out values satisfying $(h,\Pi_1)=1$ by using a simple Fundamental
Lemma sieve. For this we use the parameters $\varpi$ and $z_2$ given
by (\ref{vpd}) and (\ref{vpd1}). We then define $\xi_0(h)=\xi(h)$ if
$h<x^{1/4}$ and
\beql{Fex}
\xi_0(h)=\twosum{d\mid(h,\Pi_1)}{d<z_2}\mu(d)
\eeq
otherwise. Our immediate goal will then be the following result.
\begin{lemma}\label{NL3}
  Let
  \[\Delta(n)=\Delta(n,m)=\chi_{\cA(m)}(n)-\dz\dw^{-1}\chi_{\cB(m)}(n),\]
where $\chi_{\cA(m)}$ is the characteristic function of $\cA(m)$, for
example.  Then 
\begin{eqnarray}\label{S2bd}
 \lefteqn{S(\cA^{(k)},(4x)^{1/4})-
   \dz\dw^{-1}S(\cB^{(k)},(4x)^{1/4})}\nonumber\\
 &=& \sum_{q_1,\ldots,q_k\mid \Pi_2}\mu(q_1)\ldots\mu(q_k)\sum_{h_1,\ldots,h_k}
 \xi_0(h_1)\ldots\xi_0(h_k)\Delta(q_1h_1\ldots q_kh_k)\nonumber\\
 &&\hspace{1cm}\mbox{}+O(\dz x(\log x)^{-1}(\log\log x)^{-3/2}),
\end{eqnarray}
with suitably restricted sums over the $q_i$ and $h_i$.  Specifically
the $q_i$ are composed of
prime factors which run over disjoint dyadic intervals, and the $h_i$
also run over dyadic intervals; and no product from a subset of these dyadic
intervals falls in
any of the ranges $[x^{1/\ell-\eta},x^{1/\ell+\eta}]$ for
$4\le\ell\le\nu+2$.
\end{lemma}

The claim in the lemma is that replacing
\[\twosum{h_1,\ldots,h_k}{qh_1\ldots h_k\in\cA}\xi(h_1)\ldots\xi(h_k)\]
by
\[\twosum{h_1,\ldots,h_k}{qh_1\ldots h_k\in\cA}\xi_0(h_1)\ldots\xi_0(h_k)\]
produces a negligible error, and similarly for $\cB$. The key result which
facilitates this
is the following, which is an immediate deduction from the author's work
\cite[Lemma 15]{hbg}.
\begin{lemma}\label{hbgl}
  We have
  \[|\xi(h)-\xi_0(h)|\le\left|\twosum{d\mid(h,\Pi_1)}{d\ge z_2}\mu(d)\right|\le
  \twosum{d\mid(h,\Pi_1)}{z_2\le d<z_1 z_2}1\le \tau(h).\]
\end{lemma}

Using the lemma we then see that
\[|\xi(h_1)\ldots\xi(h_k)-\xi_0(h_1)\ldots\xi_0(h_k)|\ll
\max_j \left(|\xi(h_j)-\xi_0(h_j)|\prod_{i;\, i\not=j}\tau(h_i)\right).\]
Then, writing $h=h_1\ldots h_k$, we have
\begin{eqnarray*}
  |\xi(h_1)\ldots\xi(h_k)-\xi_0(h_1)\ldots\xi_0(h_k)|&\ll&
  \tau(h)^2  \max_j |\xi(h_j)-\xi_0(h_j)|\\
  &\ll& \tau(h)^2  \max_j\twosum{d\mid(h_j,\Pi_1)}{z_2\le d<z_1 z_2}1\\
  &\ll & \tau(h)^2 \twosum{d\mid(h,\Pi_1)}{z_2\le d<z_1 z_2}1.
\end{eqnarray*}
Since $\tau_3\mu^2*\tau_3\tau^2\le\tau^5$ the error we have to control is 
then
\[\ll\sum_{q\mid \Pi_2}\tau_k(q)\sum_{h:\, qh\in\cA}\tau_k(h)\tau(h)^2
\twosum{d\mid(h,\Pi_1)}{z_2\le d<z_1 z_2}1\ll
\twosum{d\mid \Pi_1}{z_2\le d<z_1z_2}\twosum{n\in\cA}{d|n}\tau(n)^5.\]

We have
$z_1 z_2\le x^{1/4}$ say, by (\ref{z1nu}) and (\ref{vpd1}).
Moreover, if $d\mid n$, then
$\tau(n)\le\tau(d)\tau(n/d)$, whence Lemma \ref{lsh} with $z=1$ yields
\[\twosum{n\in\cA}{d\mid n}\tau(n)^5\ll\tau(d)^5\frac{\dz x}{d}(\log x)^{31}.\]
It follows that the overall error on replacing $\xi(h)$ by $\xi_0(h)$
for $\cA(m)$ is
\[\ll \dz x(\log x)^{31}\twosum{d\mid \Pi_1}{z_2\le d<z_1 z_2}\tau(d)^5/d.\]
We bound this last sum using Rankin's trick, as follows.  For any
$\theta>0$ we have
\begin{eqnarray*}
  \twosum{d\mid \Pi_1}{z_2\le d<z_1 z_2}\tau(d)^5/d&\le&
  z_2^{-\theta}\twosum{d\mid \Pi_1}{z_2\le d<z_1 z_2}\tau(d)^5d^{-1+\theta}\\
  &\le& z_2^{-\theta}\sum_{\substack{d=1\\ d\mid \Pi_1}}^{\infty}
  \tau(d)^5d^{-1+\theta}\\
  &=& z_2^{-\theta}\prod_{p< z_1}\left(1+32 p^{-1+\theta}\right)\\
  &\le&z_2^{-\theta}\exp\left\{32\sum_{p<z_1}p^{-1+\theta}\right\}.
\end{eqnarray*}
We choose $\theta=1/\log z_1$, so that
\[\sum_{p<z_1}p^{-1+\theta}\le 3\log\log z_1,\]
(for $x$ large enough) and hence
\[\twosum{d\mid \Pi_1}{z_2\le d<z_1 z_2}\tau(d)^5/d\ll
z_2^{-1/\log z_1}(\log x)^{96}=e^{-\varpi}(\log x)^{96},\]
by (\ref{vpd1}).  Thus the error induced by replacing $\xi(h)$ with
$\xi_0(h)$ is $O(\dz x/(\log x)^2)$, say, if $\varpi$ satisfies
(\ref{vpd}).  As in the previous section,
although we have presented the argument as it applies to $\cA$,
it applies in the same way for $\cB$, and the lemma then follows.
\medskip

There is one final step that belongs in this section.  Thus far our
adjustments have affected $\cA$ and $\cB$ separately.  However, if
$h_1$, say, is large we will show that the corresponding average of
$\Delta(q_1h_1\ldots q_kh_k)$ must be negligibly small. To be
specific, we will show the following.
\begin{lemma}\label{NL4}
Let $V$ be the largest power of 2 such that $V\le x^{5/8}$. Then terms
with $h_i>V$ contribute $O(\dz x(\log x)^{-2})$ in (\ref{S2bd}).
  \end{lemma}

For the proof we consider the union of
dyadic ranges for $h_1$ covering the interval $(V,\infty)$.
If we set $f=(q_1h_1\ldots q_kh_k)/h_1$ we will
trivially have $\Delta(q_1h_1\ldots q_kh_k)=0$ unless
$f\le 8x^{3/8}$. In this latter case we find that the overall contribution is
\[\ll \sum_{f\le 8x^{3/8}}\tau(f)^6\left|\sum_{h>V}\xi_0(h)\Delta(fh)\right|
\le \sum_{f\le 8x^{3/8}}\tau(f)^6\sum_{d<z_2}\left|\sum_{g>V/d}\Delta(fdg)\right|.\]

We now recall that
\[\cA=\Z\cap(mH,mH(1+\dz)],\;\;\; \mbox{and}\;\;\;
\cB=\Z\cap(mH,mH(1+\dw)],\]
  with $\dz\le \dw$. Thus if $Vf>mH(1+\dw)$ we have $\Delta(fdg)=0$ for all
$g>V/d$. On the other hand, if $Vf\le mH$ then
\[\sum_{g>V/d}\Delta(fdg)=\left\{\frac{mH\dz}{fd}+O(1)\right\}-
\dz\dw^{-1}\left\{\frac{mH\dw}{fd}+O(1)\right\}=O(1).\]
In the remaining case $mH<Vf\le mH(1+\dw)$, and
\[\sum_{g>V/d}|\Delta(fdg)|\le\left\{\frac{mH\dz}{fd}+O(1)\right\}+
\dz\dw^{-1}\left\{\frac{mH\dw}{fd}+O(1)\right\}\ll \dz x/fd+1.\]
We therefore see that the overall contribution from terms with $h_1>V$
is
\[\ll\sum_{f\le 8x^{3/8}}\;\sum_{d<z_2}\tau(f)^6+
\twosum{f}{Vf\in\cB}\;\sum_{d<z_2}\tau(f)^6\frac{\dz x}{fd}.\]
The first sum is $\ll z_2x^{3/8+o(1)}=x^{3/8+o(1)}$, by
(\ref{vpd1}). In the second sum we have $f\gg x/V$ and
\[\twosum{f}{Vf\in\cB}\tau(f)^6\ll \dw xV^{-1}(\log x)^{63}\]
by Lemma \ref{lsh} with $z=1$. It follows that our bound is
\[\ll x^{3/8+o(1)}+ \dz\dw x(\log x)^{64}\ll \dz x(\log x)^{-2},\]
say, by (\ref{H0ch}).  This is satisfactory for the lemma.

\section{Introducing Dirichlet Polynomials}\label{secp2}

We are now interested only in the case in which $q_1\ldots q_k$ is
square-free. Suppose that
\[q_i=p_{i,1}\ldots p_{i,t_i},\;\;\;(1\le i\le k).\]
Each prime $p_{i,j}$ runs over a corresponding dyadic interval
$I_{i,j}$, and the number of possible intervals is $O(\log x)$.  Since
$t_i\le\nu$ for $i\le k$ the total number of choices for these dyadic
intervals is at most $(C\log x)^{3+3\nu}$ for some absolute constant
$C$. In the same way the variables $h_1,\ldots,h_k$ belong to dyadic
interval $J_1,\ldots,J_k$, and there are at most $(C\log x)^3$ choices
for these intervals.  By (\ref{bsd}) and (\ref{nubnd}) there are therefore 
$O(\bso)$ choices for the entire collection of intervals.
Since we have arranged that the intervals
$I_{i,j}$ are distinct, each relevant product $q=q_1\ldots q_k$
arises exactly once as $p_{i,j}$ runs over $I_{i,j}$, and each $q$
corresponds to $\tau_k(q)$ choices for the $k$-tuple
$q_1,\ldots,q_k$.  We note that
\[\tau_k(q)\le \tau_3(q)\le 3^{\omega(q)}\le 3^{3\nu}\ll\bso.\]

It will be convenient to re-label the intervals $I_{i,j}$ as
$I_1,\dots,I_t$, where $t=t_1+\ldots+t_k$, and to replace $I_j$
by $I_j\cap[z_1,(4x)^{1/4})$.  Then
\begin{eqnarray*}
  \lefteqn{\left|\sum_{q_1,\ldots,q_k\mid \Pi_2}\mu(q_1)\ldots\mu(q_k)
 \sum_{h_1,\ldots,h_k}\xi_0(h_1)\ldots\xi_0(h_k)\Delta(q_1h_1\ldots q_kh_k)
 \right|}\\
  &\ll &\bso\sum_{I_1,\ldots,I_t}\sum_{J_1,\ldots,J_k}\left|
  \twosum{p_i\in I_i,\; h_j\in J_j}{1\le i\le t,\; 1\le j\le k}
  \xi_0(h_1)\ldots\xi_0(h_k)\Delta(p_1\ldots p_th_1\ldots h_k)\right|.
\end{eqnarray*}
We also observe that $\xi_0(h)=0$ for $1<h<z_1$ so that we may replace $J_j$ by
  $J_j\cap[z_1,x^{5/8}]$. When $h_j=1$ for some index $j$ we may omit
    $J_j$ altogether, reducing $k$ by 1. As a result we may have to
allow for the possibility that $k=0$.
    Referring to Lemma \ref{lf} and (\ref{S2bd}) we see that
    \[\dz x\ll\bso\sum_{I_1,\ldots,I_t}\;\sum_{J_1,\ldots,J_k}\left|
  \twosum{p_j\in I_j}{1\le j\le t}\;\twosum{h_j\in J_j}{1\le j\le k}
  \xi_0(h_1)\ldots\xi_0(h_k)\Delta(p_1\ldots p_th_1\ldots
  h_k,m_i)\right|\]
for $i=1,\ldots,R_1$, where we now omit collections of dyadic intervals any
subset of which contains a product from any of the ranges
$[x^{1/\ell-\eta},x^{1/\ell+\eta}]$ for $4\le\ell\le\nu+2$.

Thus there is a subset of the $m_i$, with
cardinality at least $R_1\bs^{-1}$, on which the contribution
from some specific set of intervals is large. We can therefore
conclude as follows.

\begin{lemma}\label{prev}
  Suppose all parameters are as previously defined. Then there are:-
\begin{enumerate}
\item[(i)] Integers $k=0, 1, 2$ or $3$ and
  $t\in[0,3\nu]$;
\item[(ii)] Disjoint intervals $I_j=(A_j,B_j]\subseteq[z_1,(4x)^{1/4}]$
for $1\le j\le t$, with $B_j\le 2A_j$, and intervals
$J_j=(C_j,D_j]\subseteq[z_1,x^{5/8}]$ for $1\le j\le k$
  with $D_j\le 2C_j$, with the following property. For any subsets
  $\mathcal{J}_1\subseteq\{1,\ldots,t\}$ and
  $\mathcal{J}_2\subseteq\{1,\ldots,k\}$, and for any integer in the
  range $4\le\ell\le\nu+2$,  we have either 
\[\prod_{j\in\mathcal{J}_1}A_j\prod_{j\in\mathcal{J}_2}C_j\ge
x^{1/\ell+\eta}\;\;\;\mbox{or}\;\;\;
\prod_{j\in\mathcal{J}_1}B_j\prod_{j\in\mathcal{J}_2}D_j<
x^{1/\ell-\eta};\]
and
\item[(iii)] Distinct integers $m_1,\ldots, m_{R}\in[x/H,3x/H]$;
\end{enumerate}
such that
\[\bso\left|\twosum{p_j\in I_j}{1\le j\le t}\;
\twosum{h_j\in J_j}{1\le j\le k}
  \xi_0(h_1)\ldots\xi_0(h_k)\Delta(p_1\ldots p_th_1\ldots
  h_k,m_i)\right|\gg x^{1/2}\]
for $m=m_1,\ldots,m_{R}$, and with
\[{\rm Meas}(\cl{I}(x))\ll  x^{1/2+o(1)}R.\]
\end{lemma}

Here we have used the estimate $\bs\le x^{o(1)}$ in estimating
${\rm Meas}(\cl{I}(x))$.

Clearly we must have
\[(\prod A_j)(\prod C_j)\le 4x\andd(\prod B_j)(\prod D_j)\ge x\]
in order for there to be any overlap with $\cA$ or $\cB$, and we
therefore assume henceforth that
\beql{AC}
2^{-3-3\nu}x\le(\prod A_j)(\prod C_j)\le 4x.
\eeq
We now define Dirichlet polynomials
\[P_j(s)=\sum_{p\in I_j}p^{-s},\;\;\;(1\le j\le t),\]
\[F_j(s)=\sum_{h\in J_j}\xi_0(h)h^{-s},\;\;\;(1\le j\le k),\]
and
\beql{Dp}
D(s)=P_1(s)\ldots P_t(s)F_1(s)\ldots F_k(s).
\eeq
Thus $D(s)$ has coefficients supported on $[2^{-3-3\nu}x,2^{5+3\nu}x]$.

At this point it is convenient to establish a general result on the
coefficients of products of these Dirichlet polynomials.
\begin{lemma}\label{coeffest}
The coefficients of any product of distinct factors $P_j(s)$ take
values 0 and 1 only. The coefficients $c_n$ of any sub-product of
$D(s)$ satisfy $|c_n|\le\tau_7(n)$. Provided that one excludes
factors $F_j$ for which $C_j\ge x^{1/4}$, any product of powers of the
polynomials $P_j$ and $F_j$ will have coefficients $c_n$ satisfying
$|c_n|\le\nu^{\nu}\ll\bso$ for $n\le x$.
\end{lemma}

The first assertion follows from the fact that the polynomials $P_j$
have coefficients supported on primes in disjoint intervals $I_j$.
For the second claim we observe that the coefficients of $F_j$ have
size at most $\tau(n)$, so that the sub-product in question will have
coefficients dominated by those of
$\zeta(s)\left(\zeta(s)^2\right)^k$, giving us the required bound
$\tau_7(n)$. For the final assertion, we observe that if
our product of Dirichlet polynomials contains a term with $n\le x$
it can have at most $h=[\nu]$ factors, since $A_j$ and $C_j$ are at least
$z_1$.  It follows that $|c_n|\le\tau_h(n)$. Moreover the $c_n$ are
supported on products of primes $p\ge z_1$, since we are excluding the
case $C_j\ge x^{1/4}$. We then have $\Omega(n)\le\nu$, in light of the
assumption that $n\le x$.  However $\tau_h(n)$ is at most the number
of ways that a set of $\Omega(n)$
primes (distinct or not) can be partitioned into $h$ subsets, whence
$\tau_h(n)\le h^{\Omega(n)}\le \nu^{\nu}$. This completes the proof of
the lemma.

We are now ready to state the main result of this section.
\begin{Prop}\label{P2}
Let
\beql{Tdef}
T=x^{1/2}\bs^2
\eeq
and
\beql{T0def}
 T_0=\exp(\tfrac13\sqrt{\log x}).
 \eeq
 Then, in the situation of Lemma \ref{prev}, there are
complex coefficients $\zeta_j$ of modulus 1 for which the function
\[M(s)=\sum_{j=1}^R \zeta_j m_j^{-s}\]
satisfies
\beql{MB}
Rx\ll \bso \int_{T_1}^{2T_1}|D(it)M(it)|dt
\eeq
for some $T_1\in[T_0,T]$.
\end{Prop}
As explained in connection with the definition (\ref{bsd}) we think
of this final bound as involving a loss of a factor $\bso$. The
integral on the right should suggest the use of Proposition \ref{P1},
although much work must be done first.  However we
observe at this point that
the integers $m_j$ satisfy $0<m_j\le 3x/H\le T$, as required for
Proposition \ref{P1}, by virtue of (\ref{Hch}) and (\ref{Tdef}).

We begin the proof of Proposition \ref{P2}
by following the usual analysis of Perron's formula, as in Titchmarsh
\cite[Sections 3.12 and 3.19]{titch} for example. This produces
\begin{eqnarray}\label{sa}
\lefteqn{\sum_{p_j\in I_j,\, (1\le j\le t)}\;\;
\twosum{h_j\in J_j\,(1\le j\le k)}{p_1\ldots p_th_1\ldots h_k\in\cA(m)}
  \xi_0(h_1)\ldots\xi_0(h_k)}\nonumber\\
&=&
\frac{1}{2\pi i}\int_{-iT}^{iT}D(s)\frac{(1+\dz)^s-1}{s}(Hm)^sds+O(E),
\end{eqnarray}
where the error $E$ is given by
\[\sum_{2^{-3-3\nu}x\le n\le 2^{5+3\nu}x}\tau_7(n)\min\left\{
\frac{T^{-1}}{|\log(mH(1+\dz)/n)|}+
\frac{T^{-1}}{|\log mH/n|}\,,\,\log T\right\}.\]
Terms with $n<x/2$ or $n>5x$ contribute $O(\tau_7(n)/T)$ each, and
hence produce $O(2^{3\nu}x(\log x)^{6}/T)$ in total. When
$J<|mH(1+\dz)-n|\le 2J$ with
$x^{1/4}\le J\ll x$ we have
\[\frac{T^{-1}}{|\log(mH(1+\dz)/n)|}\ll (JT)^{-1}x\]
and the corresponding terms contribute $O((JT)^{-1}x(J\log^{6} x))$, by Lemma
\ref{lsh}, taking $z=1$. Summing over dyadic ranges for $J$ produces
$O(x(\log x)^7/T)$, and similarly for the contribution from
$T^{-1}/|\log mH/n|$. Finally, for terms with $|mH(1+\dz)-n|\le x^{1/4}$ or
$|mH-n|\le x^{1/4}$ we bound the minimum in $E$ by $\log T$, obtaining
a contribution $O(x^{1/4}(\log x)^{6}(\log T))$, by a further
application of Lemma \ref{lsh}. It therefore follows that
\[E\ll 2^{3\nu}x(\log x)^{6}/T+x(\log x)^{7}/T+x^{1/4}(\log x)^{6}(\log
T).\]
Our choice (\ref{Tdef}) ensures that $E\ll x^{1/2}\bs^{-1}$,
by (\ref{bsd}) and  (\ref{80}).

A similar analysis applies to $\cB(m)$, leading to the estimate
\begin{eqnarray*}
\lefteqn{\sum_{p_j\in I_j,\, (1\le j\le t)}\;\;
\twosum{h_j\in J_j\,(1\le j\le k)}{p_1\ldots p_th_1\ldots h_k\in\cB(m)}
  \xi_0(h_1)\ldots\xi_0(h_k)}\\
&=&
  \frac{1}{2\pi i}\int_{-iT}^{iT}D(s)\frac{(1+\dw)^s-1}{s}(Hm)^sds
  +O(x^{1/2}\bs^{-1}).
\end{eqnarray*}
It then follows that
\[\sum_{p_j\in I_j,\, (1\le j\le t)}
\sum_{h_j\in J_j\,(1\le j\le k)}
  \xi_0(h_1)\ldots\xi_0(h_k)\Delta(p_1\ldots p_th_1\ldots
  h_k,m)\]
  \[=\frac{1}{2\pi i}\int_{-iT}^{iT}D(s)G(s)m^sds
  +O(x^{1/2}\bs^{-1}),\]
  with
  \[G(s)=\left(\frac{(1+\dz)^s-1}{s}-\dz\dw^{-1}\frac{(1+\dw)^s-1}{s}
  \right)H^s.\]
  Now if $0\le\mu\le 1$ and $t$ is real, we have
  \[\frac{(1+\mu)^{it}-1}{it}=\mu+
  \int_0^{\mu}\int_0^{\nu}(1+\lambda)^{it-2}(it-1)d\lambda d\nu=
  \mu+O(\mu^2(1+|t|)),\]
whence  $G(it)\ll\dz\dw(1+|t|)$.  Moreover $|P_j(it)|\le A_j$ and $F_j(it)\ll
 C_j(\log x)$. We therefore deduce from (\ref{AC}) that
 \[\int_{-iT_0}^{iT_0}|D(it)G(it)|dt\ll \dz\dw x(\log x)^3T_0^2.\]
 The choice  (\ref{T0def}) shows that
 the above bound is $O(x^{1/2}\bs^{-1})$, by (\ref{dzch}),
 (\ref{H0ch}), and (\ref{81}). This allows us to conclude from part (iii) of
 Lemma \ref{prev} that
\[\bso\left|\int_{T_0<|t|\le T}D(it)G(it)m_j^{it}dt\right|\gg x^{1/2}\]
for $m_j=m_1,\ldots,m_R$, when $x$ is large enough. 

We have now reached an important stage in the argument. By choosing
suitable complex coefficients $\zeta_j$ of modulus 1 we  can write
\[\overline{\zeta_j}\int_{T_0<|t|\le T}D(it)G(it)m_j^{it}dt=
\left|\int_{T_0<|t|\le T}D(it)G(it)m_j^{it}dt\right|\]
for $1\le j\le R$, whence
\[\bso\int_{T_0<|t|\le T}D(it)G(it)\overline{M(it)}dt\gg Rx^{1/2},\]
with $M(s)$ as in Proposition \ref{P2}.  Indeed, since
\[|G(it)|=\left|\int_1^{1+\dz}v^{-it-1}dv-
\dz\dw^{-1}\int_1^{1+\dw}v^{-it-1}dv\right|\le 2\dz,\]
we have
\begin{eqnarray*}
Rx^{1/2}&\ll&\dz\bso\int_{T_0<|t|\le T}|D(it)\overline{M(it)}|dt\\
&\ll&
x^{-1/2}\bso\int_{T_0<|t|\le T}|D(it)M(it)|dt
\end{eqnarray*}
for large enough $x$, by (\ref{dzch}).  Moreover, by
dyadic subdivision there will be a value $T_1\in[T_0,T]$ such that
\begin{eqnarray*}
Rx^{1/2}&\ll&(\log x)x^{-1/2}\bso\int_{T_1\le|t|\le 2T_1}|D(it)M(it)|dt\\
&\ll& x^{-1/2}\bso\int_{T_1\le |t|\le 2T_1}|D(it)M(it)|dt.
\end{eqnarray*}
The contribution from negative $t$ has the same shape as that for
positive $t$, but with $\zeta_j$ replaced by its conjugate, so that it
suffices to consider $T_1\le t\le 2T_1$.
The proposition then follows.

It is the introduction of the
Dirichlet polynomial $M(s)$, and the estimation of mean-values
involving it, via Proposition \ref{P1}, which are the most significant
features of this paper.

\section{Extremely Large Values of Dirichlet Polynomials}

The next stage in the argument is to show that $P_j(it)$ and $F_j(it)$
cannot be extremely large.
\begin{lemma}\label{EL}
  We have
  \beql{Pub}
|P_j(it)|\le A_j\exp(-(\log x)^{1/5}),
\eeq
if $x$ is large enough.  Similarly, we have
  \beql{Fub}
|F_j(it)|\le C_j\exp(-(\log x)^{1/5}),
\eeq
if $x$ is large enough.  
\end{lemma}

For $P_j(s)$ this follows by the argument used for Lemma
19 of Heath-Brown \cite{hbg}, which handled Dirichlet polynomials
evaluated at $\tfrac12+it$ rather than $it$. Since $A_j\ge z_1$
the argument shows that
\[P_j(it)\ll A_j(z_1^{-\beta(T_1)}+T_1^{-1})(\log x)^2\]
with $\beta(T_1)$ of order $(\log T_1)^{-2/3}(\log\log T_1)^{-1/3}$,
so that $\beta(T_1)\ge(\log x)^{-3/4}$ for large $x$.  Thus
(\ref{z1nu}) and (\ref{T0def}) yield
\begin{eqnarray*}
  P_j(it)&\ll& A_j(z_1^{-(\log x)^{-3/4}}+T_0^{-1})(\log x)^2\\
  &\ll& A_j
  \{\exp(-(\log x)^{1/4}/\nu)+\exp(-\tfrac13(\log x)^{1/2})\}(\log x)^2,
\end{eqnarray*}
and (\ref{Pub}) follows by (\ref{nubnd}), if $x$ is large enough.

When $C_j\le x^{1/4}$ we have
\[F_j(it)=\sum_{C_j<n\le D_j}\xi(n)n^{-it}.\]
We write $n$ as $pm$ where $p=P^+(n)$ is the largest prime factor of $n$.
This allows us to classify terms according to the value of $m$, giving
\[|F_j(it)|\ll \sum_m
\left|\twosum{C_j/m<p\le D_j/m}{p\ge \max(P^+(m),z_1)} p^{-it}\right|.\]
The inner sum is empty unless $D_j/m\ge z_1$, and then the previous
argument shows that
\[\twosum{C_j/m<p\le D_j/m}{p>\max(P^+(m),z_1)} p^{-it}\ll
  (C_j/m)\exp(-2(\log x)^{1/5}),\]
say. The required estimate then follows on summing over $m$.

In the remaining range $C_j\ge x^{1/4}$ we have
\beql{Fex1}
F_j(it)=\sum_{C_j<n\le D_j}\xi_0(n)n^{-it}\ll
\sum_{d\le z_2}\left|\sum_{C_j/d<n\le D_j/d}n^{-it}\right|,
\eeq
by (\ref{Fex}).  When $T_1\le t\le 2T_1$ the inner sum is
\[\ll (C_j/d)^{1/2}T_1^{1/6}+(C_j/d)T_1^{-1/6}\]
by the van der Corput third derivative estimate (see Titchmarsh
\cite[Theorem 5.11]{titch}). It follows that
\[F_j(it)\ll z_2^{1/2}C_j^{1/2}T_1^{1/6}+C_jT_1^{-1/6}\log x
\ll z_2^{1/2}C_jx^{-1/8}T_1^{1/6}+C_jT_1^{-1/6}\log x,\]
since $C_j\ge x^{1/4}$.
However $T_1\le T\ll x^{1/2+o(1)}$ by (\ref{Tdef}), whence (\ref{vpd1}) yields
\[z_2^{1/2}C_jx^{-1/8}T_1^{1/6}\ll C_jx^{-1/24+o(1)};\]
and $T_1\ge T_0$, whence
\[C_jT_1^{-1/6}\log x\ll C_j\exp\{-\tfrac14\sqrt{\log x}\},\]
say, by (\ref{T0def}). The bound required for the lemma then follows.
\bigskip

We have already shown that we can take $C_j\le x^{5/8}$, and we now
reduce this bound further.
\begin{lemma}\label{8l1}
  If $C_j\ge \max(x^{1/4},T_1z_2)$ then $R=0$.
\end{lemma}
Thus we will assume henceforth that $C_j\le\max(x^{1/4},T_1z_2)$.

For the proof we apply (\ref{Fex1}). According to Titchmarsh
\cite[Theorem 4.11]{titch} we have
\[\sum_{N<n\le M}n^{-1/2-it}\ll M^{1/2}/|t|\]
uniformly for $M\ge N\ge |t|/2$, say.  In our situation we
have
\[C_j/d\ge C_j/z_2\ge T_1\ge |t|/2,\]
and it follows by partial summation that
\[\sum_{C_j/d<n\le D_j/d}n^{-it}\ll C_jd^{-1}T_1^{-1}.\]
We therefore conclude that $F_j(it)\ll C_j(\log x)/T_1$ uniformly for
$t\in[T_1,2T_1]$.  Since $C_j\le x^{5/8}$ the product $D(it)$ must
contain at least one other factor apart from $F_j(it)$.  We therefore
see from Lemma \ref{EL} that
\begin{eqnarray*}
D(it)&\ll& T_1^{-1}(\log x)\exp\{-(\log x)^{1/5}\}\prod_j A_j\prod_j C_j\\
&\ll& T_1^{-1}(\log x)\exp\{-(\log x)^{1/5}\}x,
\end{eqnarray*}
by (\ref{AC}). Since $M(it)\ll R$ we then deduce from Proposition \ref{P2}
that
\[Rx\ll\bso (\log x)\exp\{-(\log x)^{1/5}\}Rx.\]
This then shows that we must have $R=0$, in view of (\ref{81}).

\section{The Fourth Moment of $F_j(it)$}

Our next goal is the following estimate.
\begin{lemma}\label{lF4}
  Suppose that $C_j\ge x^{1/4}$.
  Let $t_1,\ldots,t_M\in[T_1,2T_1]$, and assume that $|t_m-t_n|\ge 1$ for
  $m\not=n$.  Then
  \[\sum_{m=1}^M|F_j(it_m)|^4\ll\bso T_1C_j^2.\]
  Moreover
  \[\int_{T_1}^{2T_1}|F_j(it)|^4dt\ll \bso T_1C_j^2.\]
\end{lemma}

The second claim clearly follows from the first.  The lemma
would still be true when $C_j<x^{1/4}$, but we only need the
lemma for $C_j\ge x^{1/4}$. It would be quite easy to establish an
estimate of the above form with an additional factor $z_2^4$, say,
using the classical fourth moment estimate for the Riemann
zeta-function; but unfortunately $z_2\not=O(\bso)$. The key input for
the proof is therefore the following estimate, which is an immediate
corollary of Theorem 1 of Bettin, Chandee and Radziwi\l\l\, \cite{BCR}.
\begin{lemma}\label{bcr}
  Let
  \[A(s)=\sum_{n\le N}a_nn^{-s}\]
  with $|a_n|\le\tau_3(n)$ and $N\le T^{1/2+1/67}$.  Then
  \[\int_T^{2T}|\zeta(\tfrac12+it)A(\tfrac12+it)|^2dt\ll
  T\log T\sum_{m,n}\frac{|a_ma_n|}{[m,n]}+T.\]
\end{lemma}
Indeed
\begin{eqnarray*}
  \sum_{m,n\le N}\frac{\tau_3(m)\tau_3(n)}{[m,n]}&\le&\sum_{d\le N}
  \twosum{m,n\le N}{d\mid m,n}\frac{\tau_3(m)\tau_3(n)}{mn/d}\\
  &\le &\sum_{d\le N}\frac{\tau_3(d)^2}{d}
  \left\{\sum_{u\le N/d}\frac{\tau_3(u)}{u}\right\}^2\\
  &\ll &\sum_{d\le N}\frac{\tau_3(d)^2}{d}
  \left\{(\log N)^3\right\}^2\\
  &\ll &(\log N)^{15},
\end{eqnarray*}
whence
\beql{bcre}
\int_T^{2T}|\zeta(\tfrac12+it)A(\tfrac12+it)|^2dt\ll T(\log T)^{16}.
\eeq
\bigskip

We begin our proof of Lemma \ref{lF4} by writing
\[Z(s)=\twosum{d<z_2}{d\mid\Pi_1}\mu(d)d^{-s}\;\;\;
\mbox{and}\;\;\; F(s)=Z(s)\zeta(s).\]
It then follows from (\ref{Fex}) that the coefficients of $F_j(s)$ and
$F(s)$ agree for $C_j<n\le D_j$. Using Perron's formula (Titchmarsh
\cite[Sections 3.12 and 3.19]{titch} for example) we then deduce that
\[F_j(it)=\frac{1}{2\pi i}\int_{5/4-iT_1/2}^{5/4+iT_1/2}F(s+it)
\frac{D_j^s-C_j^s}{s}ds+O(E(C_j))+O(E(D_j))\]
for $T_1\le t\le 2T_1$, with
\[E(A)=\sum_{n=1}^\infty\tau(n)(A/n)^{5/4}\min\left\{
\frac{T_1^{-1}}{|\log(A/n)|}\,,\,\log T_1\right\}.\]
Since $C_j\le\max(x^{1/4},T_1z_2)$ we have $C_j\le T_1z_2$, so that
$T_1\ge C_j^{7/8}$, say. An analysis similar to that used for
(\ref{sa}) then shows that the error terms satisfy
\[E(C_j)+E(D_j)\ll C_j^{1/2}.\]

We proceed to move the line of integration to Re$(s)=\tfrac12$,
incurring an error
\[\left(\frac{C_j^{5/4}}{T_1}+\frac{C_j^{1/2}T_1^{1/4}}{T_1}\right)z_2\ll
C_j^{1/2},\]
say, whence
\[F_j(it)\ll C_j^{1/2}\left(1+\int_{-T_1/2}^{T_1/2}
|F(\tfrac12+i(\tau+t))|\frac{d\tau}{1+|\tau|}\right).\]
We then find via H\"{o}lder's inequality that
  \[\sum_{m=1}^M|F_j(it_m)|^4\ll
  C_j^2\left\{T_1+(\log T_1)^3\sum_{m=1}^M\int_{-T_1/2}^{T_1/2}
  |F(\tfrac12+i(\tau+t_m))|^4
  \frac{d\tau}{1+|\tau|}\right\}.\]
However 
\begin{eqnarray*}
\int_{-T_1/2}^{T_1/2}|F(\tfrac12+i(\tau+t_m))|^4\frac{d\tau}{1+|\tau|}
&=&
\int_{t_m-T_1/2}^{t_m+T_1/2}|F(\tfrac12+i\tau)|^4
\frac{d\tau}{1+|\tau-t_m|}\\
&\le &\int_{T_1/2}^{5T_1/2}|F(\tfrac12+i\tau)|^4
\frac{d\tau}{1+|\tau-t_m|},
\end{eqnarray*}
and
\[\sum_{m=1}^M\frac{1}{1+|\tau-t_m|}\ll \log T_1,\]
in view of the spacing condition on the points $t_m$.  We therefore
conclude that
\[\sum_{m=1}^M|F_j(it_m)|^4\ll C_j^2\left\{T_1+(\log T_1)^4\int_{T_1/2}^{5T_1/2}
  |F(\tfrac12+i\tau)|^4d\tau\right\}.\]
  \medskip

  Our task is now to estimate the fourth power moment of the function
  $F(\tfrac12+i\tau)=Z(\tfrac12+i\tau)\zeta(\tfrac12+i\tau)$,
  using Lemma \ref{bcr}. For
  this we will employ the approximate functional equation for $\zeta(s)$
  as given by Titchmarsh \cite[Theorem 4.13]{titch}.  This yields
  \[\zeta(\tfrac12+i\tau)\ll |Z_1(\tfrac12+i\tau)|+
  |Z_2(\tfrac12+i\tau;\tau)|+1,\]
  with
  \beql{Xch}
  Z_1(s)=\sum_{n\le X}n^{-s}\;\;\; (X=T_1^{1/2+1/68})
  \eeq
  and
  \[Z_2(s;\tau)=\sum_{n\le 2\pi\tau/X}n^{-s}.\]
  We apply this to  two factors $\zeta(s)$, whence
  \beql{123}
  \sum_{m=1}^M|F_j(it_m)|^4\ll C_j^2\{T_1+(\log T_1)^4(I_1+I_2+I_3)\},
  \eeq
  with
  \[I_1=\int_{T_1/2}^{5T_1/2}|\zeta(\tfrac12+i\tau)Z_1(\tfrac12+i\tau)
  Z(\tfrac12+i\tau)^2|^2d\tau,\]
  \[I_2=\int_{T_1/2}^{5T_1/2}|\zeta(\tfrac12+i\tau)Z_2(\tfrac12+i\tau;\tau)
  Z(\tfrac12+i\tau)^2|^2d\tau,\]
  and
\[I_3=\int_{T_1/2}^{5T_1/2}|\zeta(\tfrac12+i\tau)
  Z(\tfrac12+i\tau)^2|^2d\tau.\]
The first and third of these can be handled immediately by
(\ref{bcre}), giving bounds $O(\bso T_1)$.

For $I_2$ we use Cauchy's inequality together with the usual 
fourth moment estimate for the Riemann zeta-function
(Titchmarsh \cite[(7.6.1)]{titch} for example) to show that
\beql{I24}
I_2^2\le\left\{\int_{T_1/2}^{5T_1/2}|\zeta(\tfrac12+i\tau)|^4d\tau\right\}
I_4\ll T_1(\log T_1)^4I_4,
\eeq
with
\[I_4=\int_{T_1/2}^{5T_1/2}
|Z_2(\tfrac12+i\tau;\tau)Z(\tfrac12+i\tau)^2|^4d\tau.\]
If we expand $Z_2^2Z^4\overline{Z_2}^2\overline{Z}^4$ we get a sum
over 12-tuples $(m_1,\ldots,m_4,n_1,\ldots,n_8)$ in which
$m_j\le 2\pi\tau/X$ and $n_j<z_2$. We then have to examine
\[\int_{\tau\in[T_1/2,5T_1/2],\,\tau\ge(X/2\pi)\max(m_j)}
\left(\frac{U}{V}\right)^{i\tau}d\tau,\]
where $U=m_1m_2n_1n_2n_3n_4$ and $V=m_3m_4n_5n_6n_7n_8$.  Since the
range for $\tau$ is a sub-interval of $[T_1/2,5T_1/2]$ the integral is
of order $|\log U/V|^{-1}$ whenever $U\not=V$, and is of order $T_1$
otherwise. Each value of $U$ occurs at most $\tau_6(U)$ times, and
similarly for $V$.  Moreover since $\tau\le 5T_1/2$ we have
$m_i\ll T_1/X$, whence $U,V\ll (T_1/X)^2z_2^4\ll T_1$ via our choice
(\ref{Xch}) for $X$. It follows that
\[I_4\ll T_1\sum_{U\ll T_1}\frac{\tau_6(U)^2}{U}
+\sum_{U\not=V\ll T_1}\frac{\tau_6(U)\tau_6(V)}{(UV)^{1/2}}|\log U/V|^{-1}.\]
The first sum is $O((\log T_1)^{36})$.  For the second we note that
$2ab\le a^2+b^2$ for real $a,b$ whence
\begin{eqnarray*}
\lefteqn{\sum_{U\not=V\ll T_1}
    \frac{\tau_6(U)\tau_6(V)}{(UV)^{1/2}}|\log U/V|^{-1}}\hspace{3cm}\\
  &\le&\tfrac12
  \sum_{U\not=V\ll T_1}\left(\frac{\tau_6(U)^2}{U}
  +\frac{\tau_6(V)^2}{V}\right)|\log U/V|^{-1}\\
  &=& \sum_{U\not=V\ll T_1}\frac{\tau_6(U)^2}{U}|\log U/V|^{-1},
  \end{eqnarray*}
by symmetry.  However
\[\twosum{V\ll T_1}{V\not=U}|\log U/V|^{-1}\ll T_1\log T_1,\]
whence $I_4\ll T_1(\log T_1)^{37}$. Thus $I_2\ll\bso T_1$ by
(\ref{I24}), so that (\ref{123}) yields the required estimate for Lemma
\ref{lF4}.

\section{Large Values of Dirichlet Polynomials}

In this section we handle moderately large values of $P_j(it)$ and
$F_j(it)$. For this section only, it will be convenient to define
\[P_{t+j}(s)=F_j(s)\;\; \mbox{and}\;\; A_{t+j}=C_j,\;
B_{j+t}=D_j\;\; \mbox{for}\;\;\ 1\le j\le k. \]
We proceed by covering the range $[T_1,2T_1]$ with unit intervals
$[n,n+1]$, and
examining the contribution from those where one has 
\beql{45}
\sup_{[n,n+1]}|P_j(it)|> A_j^{4/5}\;\;\;\mbox{for some}\; j\le t+k.
\eeq
Consider the set $\cl{N}$ of integers $n$ with
$[T_1]\le n\le [2T_1]$ for which (\ref{45}) holds, and such that
\[V_0<\sup_{[n,n+1]}|D(it)|\le 2V_0\]
for a given $V_0\ge 1$, where $D(s)$ is given by (\ref{Dp}).
Our goal is to prove the following bound.
\begin{lemma}\label{clNbnd}
  We have
  \[\card\cl{N}\ll xV_0^{-1}\exp\{-(\log x)^{1/6}\}.\]
\end{lemma}

Before proceeding to the proof of the lemma we note the following
consequence. By using a
dyadic subdivision into values of $V_0$ which are powers of 2, we will
have
\[\int |D(it)|dt\ll x(\log x)\exp\{-(\log x)^{1/6}\}\]
where the integral is over relevant intervals $[n,n+1]$ such that
(\ref{45}) holds. (Note that intervals where $V_0\le 1$
contribute a total $O(T_1)$, which is satisfactory.)
We may compare this with the lower
bound in Proposition \ref{P2}.  Since $|M(it)|\le R$ for all $t$, and 
\[x(\log x)\exp\{-(\log x)^{1/6}\}\ll x\bs^{-1},\]
by (\ref{81}), we see that intervals $[n,n+1]$ where
(\ref{45}) holds make a negligible contribution in
Proposition \ref{P2}.  Thus,
in subsequent work we will be able to assume  that for all
relevant $t$ we have $|P_j(it)|\le A_j^{4/5}$ and
$|F_j(it)|\le C_j^{4/5}$.
\medskip

For the proof of the lemma we write
\[x_0=\prod_1^{t+k} A_j\le 4x,\]
by (\ref{AC}).  For each $j$ we choose $V_j=V_j(n)$ to be a power of 2 with
\[V_j/2<\sup_{[n,n+1]}|P_j(it)|\le V_j,\]
and we choose $j=j_0=j_0(n)$ such that
\[\sigma=\frac{\log V_j}{\log A_j}\]
is maximal. In particular we have $\sigma> 4/5$ by (\ref{45}).
Moreover 
\[A_j^{1-\sigma}\ge\tfrac12\exp((\log x)^{1/5})\]
by Lemma \ref{EL}. Thus 
\[x^{1-\sigma}\ge\exp((\log x)^{1/5}), \]
since $A_j\le x^{5/8}$ by part (ii) of Lemma \ref{prev}.  It follows that
\beql{sub}
1-\sigma\ge (\log x)^{-4/5}.
\eeq
We also have
\[V_0\le\sup_{[n,n+1]}|D(it)|\le\prod_j\sup_{[n,n+1]}|P_j(it)|\le
\prod_j V_j\le \left(\prod_j A_j\right)^{\sigma}=x_0^{\sigma},\]
so that
\beql{V0}
x^{\sigma}\gg x_0^{\sigma}\ge V_0.
\eeq

We subdivide the integers $n\in\cl{N}$ according to the
values of $j_0(n)$ and $V_{j_0}$, producing $O((\log x)^2)$
subsets. There is thus a choice of $j_0$ and $V_{j_0}$ for which the
corresponding subset, $\cl{N}_1$ say, satisfies
$\card\cl{N}\ll(\log x)^2\card\cl{N}_1$, and such that
\[\sup_{[n,n+1]}|P_{j_0}(it)|\ge \tfrac12 A_{j_0}^{\sigma}\]
for $n\in\cl{N}_1$. It will be typographically convenient to drop the
subscript $j_0$, and to write
\[\sup_{[n,n+1]}|P(it)|\ge \tfrac12 A^{\sigma}.\]
We then get a succession of points $t_n\in[n,n+1]$ where the suprema
are attained, and by restricting either to even $n$ in $\cl{N}_1$ or
to odd $n$, and then re-labeling, we obtain a set of points
$t_1,\ldots,t_K \in[T_1,2T_1]$, with $\card\cl{N}\ll(\log x)^2K$, such
that $|t_i-t_j|\ge 1$ for $i\not=j$, and with
\[|P(it_m)|\ge\tfrac12 A^{\sigma}\;\;\; (m=1,\ldots,K).\]

We first dispose of the case in which $P(s)=F_j(s)$ with
$A=C_j>x^{1/4}$. In this situation Lemma \ref{lF4} shows that
$KA^{4\sigma}\ll\bso TA^2$. Since $\sigma>4/5$ and
$C_j\ge x^{1/4+\eta}$ by part (ii) of Lemma \ref{prev} we deduce via
(\ref{Tdef}) that
\[\card{N}\ll\bso TA^{2-4\sigma}\ll\bs^3 x^{1/2+(2-4\sigma)(\tfrac14+\eta)}
\ll\bs^3 x^{1-\sigma-\eta}.\]
This is enough for Lemma \ref{clNbnd}, in view of (\ref{81}).

We turn now to the case in which (\ref{45}) holds for a polynomial
with $A<x^{1/4}$, so that $z_1\le A\le x^{1/4-\eta}$, by part (ii) of
Lemma \ref{prev}. We will  use the standard theory of large values estimates
for Dirichlet polynomials, considering two separate sub-cases.
Suppose firstly that
\[A\le x^{1/4-\eta};\;\mbox{and either}\; A\ge
x^{3/14}\;\mbox{or}\; \sigma\ge \frac{9}{10}.\]
We choose a non-negative integer $w$ such that
\[A^w\le x^{1/2-2\eta}< A^{w+1}.\]
Thus $w\ge 2$, since $A\le x^{1/4-\eta}$, and therefore
\beql{Arng2}
A^w\ge (x^{1/2-2\eta})^{w/(w+1)}\ge x^{1/3-2\eta}.
\eeq
Since $A\ge z_1=(4x)^{1/\nu}$ we will have $w\le\nu$. The Dirichlet
polynomial $P(s)^w$ has coefficients $c_n$ supported on integers
$n\le (2A)^w\le x$, and Lemma \ref{coeffest} shows that
$|c_n|\le \nu^{\nu}$.

We will now apply the following ``Large Values Estimate''.
\begin{lemma}\label{MLV}
  Let $t_1,\ldots,t_K\in[0,T]$ with $|t_{i}-t_j|\ge 1$ for
  $i\not=j$. Suppose we have complex coefficients $c_n$ such that
  \[\left|\sum_{n=1}^N c_n n^{-it_j}\right|\ge V\]
  for $1\le j\le K$.  Then
  \[K\ll GNV^{-2}+G^3NTV^{-6}(\log NT)^2,\]
  where
  \[G=  \sum_{n\le N}|c_n|^2.\]
\end{lemma}
This follows from Huxley \cite[(2.9)]{hux}.

In our situation we take $\sum c_n n^{-s}=P(s)^w$, with
$N=(2A)^w$, $V=2^{-w}A^{w\sigma}$, and $G\le(2A)^w\nu^{2\nu}$.  This leads to
\[K\ll \{A^{w(2-2\sigma)}+TA^{w(4-6\sigma)}\}2^{10w}\nu^{6\nu}(\log x)^2
\ll  \bso\{A^{w(2-2\sigma)}+TA^{w(4-6\sigma)}\},\]
by (\ref{80}).
It should be emphasized that this holds uniformly with respect to $w$.
If $A\ge x^{3/14}$ then since $\sigma>4/5$ we see that
(\ref{Tdef}) yields
\[A^{w(4\sigma-2)}\ge A^{6w/5}\ge A^{12/5}\ge x^{18/35}\ge T,\]
whence $A^{w(2-2\sigma)}\ge TA^{w(4-6\sigma)}$. On the other hand, if
$\sigma\ge 9/10$ then
\[A^{w(4\sigma-2)}\ge A^{8w/5}\ge x^{\tfrac85(\tfrac13-2\eta)}
\ge x^{8/15-4\eta}\ge T,\]
by (\ref{Arng2}), and again we find that $A^{w(2-2\sigma)}\ge TA^{w(4-6\sigma)}$.
Thus, under our current assumptions, we have
\[K\ll \bso A^{w(2-2\sigma)},\]
whence
\[\card\cl{N}\ll\bso x^{1-\sigma-\eta(2-2\sigma)}.\]
However $x^{1-\sigma}\ll xV_0^{-1}$ by (\ref{V0}), and (\ref{sub}) yields
\[x^{-\eta(2-2\sigma)}\ll \exp(-2\eta(\log x)^{1/5}).\]
Here $\eta(\log x)^{1/5}\ge(\log x)^{1/6}$, say, by (\ref{etanu}). Thus
(\ref{81}) yields
\[\bso x^{-\eta(2-2\sigma)}\ll\exp(-(\log x)^{1/6}),\]
  and we obtain the bound required for Lemma \ref{clNbnd} in the
  current case.

  The final situation we examine is that in which we have $A\le x^{3/14}$ and
  $4/5\le \sigma\le 9/10$.  We begin as before, but now choosing
 $w$ so that
  \[A^w\le x^{15/31}<A^{w+1}.\]
  Thus $w\ge 2$ and hence $A^w>x^{10/31}$. Instead of
  Lemma \ref{MLV} we use the following estimate.   

\begin{lemma}\label{jutl}
  Let $t_1,\ldots,t_J\in[\tau_0,\tau_0+\tau]$ with $|t_i-t_j|\ge 1$
  for
  $i\not=j$. Then for any complex coefficients $a_m$ and any fixed $\ep>0$
  we have
  \begin{eqnarray*}
    \lefteqn{\left\{\sum_{j\le J}
      \left|\sum_{M<m\le 2M}a_m m^{-it_j}\right|\right\}^2}\\
&  \ll_{\ep} &
  \tau^{\ep}(JM+J^{11/6}\tau^{1/2}+J^{23/12}\tau^{1/12}M^{1/2})
  \sum_{M<m\le 2M}|a_m|^2.
  \end{eqnarray*}
  \end{lemma}

This follows from the analysis in Section 3 of Jutlia \cite{1114}, taking
$k=3$.  We cover the range
$[T_1,2T_1]$ with $O(1+T_1/\tau)$ subintervals of length at most
$\tau$, whence some such subinterval contains $J$ points $t_j$, with
\[K\ll (1+T_1/\tau)J.\]
We proceed to split the sum $P(s)^w=\sum c_mm^{-s}$ into dyadic ranges, and 
deduce that there is some $M\le (2A)^w$ such that
\[(\tfrac12 A^{\sigma})^wJ\ll (\log x)
\sum_{j\le J}\left|\sum_{M<m\le 2M}c_m m^{-it_j}\right|.\]
Lemma  \ref{jutl} then shows that
\begin{eqnarray*}
\lefteqn{(\tfrac12 A^{\sigma})^{2w}J^2}\\
&\ll_{\ep}&(\log x)^2x^{\ep}
(J(2A)^w+J^{11/6}\tau^{1/2}+J^{23/12}\tau^{1/12}(2A)^{w/2})\sum_{n\le (2A)^w}|c_n|^2.
\end{eqnarray*}
Since $|c_n|\le \nu^{\nu}=x^{o(1)}$ by (\ref{nubnd}) this simplifes to
give
\[A^{2w\sigma}J^2\ll
x^{o(1)}(JA^w+J^{11/6}\tau^{1/2}+J^{23/12}\tau^{1/12}A^{w/2})A^w,\]
and hence
\[J\ll
x^{o(1)}(A^{w(2-2\sigma)}+\tau^3A^{w(6-12\sigma)}+\tau A^{w(18-24\sigma)}).\]
However
\[\left\{A^{w(2-2\sigma)}\right\}^{2/3}\left\{\tau^3A^{w(6-12\sigma)}\right\}^{1/3}
=\tau A^{w(10-16\sigma)/3}\ge \tau A^{w(18-24\sigma)}\]
for $\sigma\ge 4/5$, and so the final term may be dropped.
Now, since
\[\card\cl{N}\ll(\log x)^2K\ll(\log x)^2(1+T_1/\tau)J\ll
x^{o(1)}(1+x^{1/2}/\tau)J\]
we have
\[\card\cl{N}\ll x^{o(1)}\left\{1+\frac{x^{1/2}}{\tau}\right\}
(A^{w(2-2\sigma)}+\tau^3A^{w(6-12\sigma)}).\]
We choose
\[\tau=A^{w(10\sigma-4)/3},\]
whence
\[\card\cl{N}\ll
x^{o(1)}\left\{A^{w(2-2\sigma)}+x^{1/2}A^{w(10-16\sigma)/3}\right\}.\]
Since $A^w\le x^{15/31}$ and $\sigma\le 9/10$ we have
\[A^{w(2-2\sigma)}\le x^{1-\sigma-(1-\sigma)/31}\le x^{1-\sigma-1/310}.\]
Moreover, since $A^w\ge x^{10/31}$ and $\sigma\ge 4/5$ we have
\[x^{1/2}A^{w(10-16\sigma)/3}\le x^{1/2+10(10-16\sigma)/93}\le
x^{1-\sigma-1/930},\]
on noting that
\[\frac{1}{2}+\frac{10(10-16\sigma)}{93}\le 1-\sigma-\frac{1}{930}\]
on $[4/5,1]$, with equality at the lower endpoint.  These estimates
give suitable bounds for $\card\cl{N}$ in this final case. This
completes the proof of Lemma \ref{clNbnd}.

\section{Factors of Length Below $x^{1/4}$ --- The Key Proposition}

This section will be devoted to the proof of a general estimate which
will be used to handle a number of different cases.
We suppose that we have arranged the factors of $D(s)$ into three
groups, so that
\[D(s)=P_1(s)\ldots P_t(s)F_1(s)\ldots F_k(s)=A(s)B(s)C(s).\]
We suppose further that any factor $F_j(s)$ of $A(s)$ has length at
most $x^{1/4}$. We write our Dirichlet polynomials as
\[A(s)=\sum_{A<n\le 2^{3\nu+3}A}a_n n^{-s},\;\;\;
B(s)=\sum_{B<n\le 2^{3\nu+3}B}b_n n^{-s},\]
and
\[C(s)=\sum_{C<n\le 2^{3\nu+3}C}c_nn^{-s},\]
where $|a_n|\le\nu^\nu$ and $|b_n|\le\tau_7(n)\le\tau(n)^6$, by
Lemma \ref{coeffest}. We may assume that $2^{-3-3\nu}x\le ABC\le 4x$,
as in (\ref{AC}).

We now have the following result.
\begin{Prop}\label{pb}
Suppose that
  \beql{ABC}
  A\le Bx^{o(\eta)},\;\;\; B\le x^{1/2}x^{o(1)}\;\;\;
  \mbox{and}\;\;\; C^2\le Ax^{o(\eta)}.
  \eeq
    Then if $Rx\ll \bso I$ with
 \[I=\int_{T_1\le t\le 2T_1;\, |C(it)|\le C^{4/5}}|A(it)B(it)C(it)M(it)|dt\]
we have
\[R\ll_{\ep}x^{1/10+\ep},\]
for any fixed $\ep>0$.
\end{Prop}
The reader will see that the result holds under somewhat weaker but
more complicated
conditions.  However the above suffices for our
needs. Moreover one
sees that the exponent $1/10$ corresponds to the situation in which $A,B$
and $C$ are roughly $x^{2/5}, x^{2/5}$ and $x^{1/5}$.  This is the
critical case for our theorem. Clearly (\ref{ABC}) implies that
$C<x^{1/4}$, since $ABC\le 4x$.

We start by using Cauchy's inequality to show that $I\le(I_1I_2)^{1/2}$, where
\begin{eqnarray}\label{dom}
  I_1&=&\int_0^T|A(it)M(it)|^2dt\nonumber\\
  &\ll_{\ep}&
\left(A^2R^2+(AT)^{\ep}\{ART+AR^{7/4}T^{3/4}\}\right)\max_n|a_n|^2,
\end{eqnarray}
for any fixed $\ep>0$, by Proposition \ref{P1}, and
\[I_2=\int_{T_1\le t\le 2T_1;\, |C(it)|\le C^{4/5}}|B(it)C(it)|^2dt.\]
In (\ref{dom}) we have $|a_n|^2\le\nu^{2\nu}\ll \bso$.
Thus depending on which of the three terms in
(\ref{dom}) dominates we find that
\[\left\{Rx\right\}^2\ll_{\ep}\bso A^2R^2 I_2,\]
or
\[\left\{Rx \right\}^2\ll_{\ep} \bso (AT)^{\ep}ARTI_2,\]
or
\[\left\{Rx\right\}^2\ll_{\ep}\bso (AT)^{\ep}AR^{7/4}T^{3/4}I_2.\]
Rearranging these leads to
\beql{d1}
x^2\ll_{\ep} \bso A^2I_2,
\eeq
or
\beql{d2}
R\ll_{\ep}x^{-2+2\ep}ATI_2,
\eeq
or
\beql{d3}
R^{1/4}\ll_{\ep}x^{-2+2\ep}AT^{3/4}I_2.
\eeq
We will show that (\ref{d1}) cannot happen, for large $x$, and that
both (\ref{d2}) and (\ref{d3}) produce the described bound for $R$.

To estimate $I_2$ we cover the range $[T_1,2T_1]$ with
intervals $[n,n+1]$ and focus attention either on even values of $n$,
or on odd values, depending on which case makes the larger contribution.
For each such interval we choose a point $t_n$
for which $|B(it)C(it)|$ is maximal, subject to the condition that
$|C(it)|\le C^{4/5}$. Intervals in which $B(it_n)$ or $C(it_n)$ is of
order $x^{-1}$, say, contribute at most $O(T_1)$ to $I_2$. We
subdivide the remaining points further into
$O(\log^2 x)$ classes according to the dyadic ranges
\beql{dy}
U_1<|B(it_n)|\le 2U_1,\;\;\; U_2<|C(it_n)|\le 2U_2,
\eeq
in which $|B(it_n)|$ and $|C(it_n)|$ lie. After renumbering the points
$t_n$ we find that
\beql{I2b}
I_2\ll T_1+U_1^2U_2^2K\log^2 x,
\eeq
where (\ref{dy}) holds for $1\le n\le K$.

Our task now is to estimate $K$, for which we will use Huxley's large
values estimate, given by Lemma \ref{MLV}, and the mean value
estimate of Montgomery \cite[Theorem 7.3]{mont}, taking $Q=1$,
$\chi=1$, $\delta=1$.  This latter result produces the following bound.
\begin{lemma}\label{mve}
  Under the assumptions of Lemma \ref{MLV} we have
  \[K\ll G(N+T)V^{-2}\log N.\]
  \end{lemma}
We may apply Lemmas \ref{MLV} and \ref{mve} to $B(it)$, noting that
\[\sum |b_n|^2\ll\sum\tau(n)^{12}\ll (\log x)^{4095}B\ll\bso B,\]
say, by (\ref{80}), to show that
\beql{K1}
K\ll \{B^2U_1^{-2}+\min(BTU_1^{-2},B^4TU_1^{-6})\}\bso.
\eeq

We first consider the case in which the term $B^2U_1^{-2}$ in (\ref{K1})
dominates. Then (\ref{I2b}) becomes
\beql{I2bb}
I_2\ll T+B^2U_2^2\bso\ll T+B^2C^{8/5}\bso,
\eeq
since $U_2\le C^{4/5}$, whence (\ref{d1}) would produce
\[x^2\ll_{\ep}\bso A^2(T+B^2C^{8/5}).\]
However $ABC\le 4x$, and (\ref{ABC}) yields $A\le Bx^{o(\eta)}$,
whence $A\ll x^{2/3}$, for example. Using (\ref{Tdef}) we then find that
\beql{ne1}
x^2\ll_{\ep} \bso(x^{11/6}+x^2C^{-2/5})
\ll_{\ep}\bso (x^{11/6}+x^2z_1^{-2/5}),
\eeq
say.  This would provide a contradiction, by (\ref{81}).
Thus (\ref{d1}) does not hold, when
the term $B^2U_1^{-2}$ in (\ref{K1}) dominates. Moreover, if $B^2U_1^{-2}$
dominates, then the bound (\ref{I2bb})
shows that (\ref{d2}) and (\ref{d3}) reduce to
\[R\ll_{\ep}x^{-2+2\ep}AT(T+B^2C^{8/5}\bso)\ll_{\ep}
x^{-3/2+3\ep}A(x^{1/2}+B^2C^{8/5})\]
and
\[R^{1/4}\ll_{\ep}x^{-2+2\ep}AT^{3/4}(T+B^2C^{8/5}\bso)\ll_{\ep}
x^{-13/8+3\ep}A(x^{1/2}+B^2C^{8/5})\]
respectively. However our assumptions (\ref{ABC})
give $A\le x^{2/3}$ as before, and
\begin{eqnarray*}
AB^2C^{8/5}&\ll_{\ep}& AB^2A^{1/5}C^{6/5}x^{\ep}\\
&=&(ABC)^{6/5}B^{4/5}x^{\ep}\\
&\ll_{\ep}& x^{6/5}(x^{1/2+\ep})^{4/5}x^{\ep}\\
&\le&x^{8/5+2\ep},
\end{eqnarray*}
so that $R\ll_{\ep} x^{1/10+5\ep}$ in either case. This is
sufficient, on re-defining $\ep$.

For the remainder of the proof we may therefore assume that (\ref{K1})
reduces to
\beql{KC}
K\ll \min(BTU_1^{-2},B^4TU_1^{-6})\bso.
\eeq
In addition to considering mean and large values of $B(s)$ we can use
the Dirichlet polynomial
\[C(s)^w=\sum_{n\le N}c_{n}n^{-s},\]
where the integer $w\ge 2$ is chosen so that
\beql{Crng}
C^{2w-1}\le x<C^{2w+1}.
\eeq
Thus $2w-1\le\nu$, by (\ref{z1nu}), and $N\le 2^w C^w\le 2^{\nu}C^w$.
Moreover $|c_{n}|\le \nu^{\nu}$ by Lemma \ref{coeffest}.
Here we have
\[\sum|c_{n}|^2\le 2^{\nu}\nu^{2\nu}C^w,\]
whence Lemma \ref{mve} yields
\[K\ll \{C^{2w}+C^wT\}U_2^{-2w}\bso.\]

In view of part (ii) of Lemma \ref{prev} we see that the range $[C,2C]$
does not overlap any interval $[x^{1/\ell-\eta},x^{1/\ell+\eta}]$ with
$4\le \ell\le\nu+2$. Hence (\ref{Crng}) implies that we have
\beql{Crng+}
x^{1/(2w+1)+\eta}\ll C\ll x^{1/(2w-1)-\eta}.
\eeq
Moreover (\ref{ABC}) yields
\[C=(ABC)^{1/5}(C^2/A)^{2/5}(A/B)^{1/5}\ll x^{1/5+o(\eta)},\]
and it follows that $w\not=2$. 

By virtue of (\ref{KC}) we have
\[K\ll \min\left(BTU_1^{-2}\,,\,B^4TU_1^{-6}\,,\,
\{C^{2w}+C^wT\}U_2^{-2w}\right)\bso,\]
and therefore
\begin{eqnarray*}
K&\ll&\left(BTU_1^{-2}\right)^{1-3/2w}\left(B^4TU_1^{-6}\right)^{1/2w}
\left(\{C^{2w}+C^wT\}U_2^{-2w}\right)^{1/w}\bso\\
&\ll&\{T^{1-1/w}+TC^{-1}\}B^{1+1/2w}C^2U_1^{-2}U_2^{-2}\bso.
\end{eqnarray*}
Thus (\ref{I2b}) becomes
\begin{eqnarray}\label{nx}
  I_2&\ll& T+\{T^{1-1/w}+TC^{-1}\}B^{1+1/2w}C^2\bso\nonumber\\
  &\ll& \{T^{1-1/w}+TC^{-1}\}B^{1+1/2w}C^2\bso.
  \end{eqnarray}
We first use this to examine (\ref{d1}), which produces
\begin{eqnarray}\label{RR}
  x^2&\ll_{\ep}& \bso A^2\{T^{1-1/w}+TC^{-1}\}
B^{1+1/2w}C^2\nonumber\\
&\ll_{\ep}& x^{o(\eta)}A^2\{x^{1/2-1/2w}+x^{1/2}C^{-1}\}
B^{1+1/2w}C^2,
\end{eqnarray}
by (\ref{Tdef}) and (\ref{81}).
However, since $A\le Bx^{o(\eta)}$ and $ABC\le 4x$, the inequalities
(\ref{Crng+}) yield
\begin{eqnarray*}
  A^2B^{1+1/2w}C^2&\le& (AB)^{3/2+1/4w}C^2x^{o(\eta)}\\
  &=&(ABC)^{3/2+1/4w}C^{(2w-1)/4w}x^{o(\eta)}\\
  &\ll& x^{3/2+1/4w}\cdot x^{1/4w-(2w-1)\eta/4w}x^{o(\eta)},
\end{eqnarray*}
so that the overall contribution of this term to (\ref{RR}) is
\[\ll x^{o(\eta)}\cdot x^{1/2-1/2w}\cdot x^{3/2+1/2w-\eta/4}=o(x^2).\]
Similarly, we find that
\begin{eqnarray*}
  A^2B^{1+1/2w}C&\le& (AB)^{3/2+1/4w}Cx^{o(\eta)}\\
  &=&(ABC)^{3/2+1/4w}C^{-(2w+1)/4w}x^{o(\eta)}\\
  &\ll& x^{3/2+1/4w}\cdot x^{-1/4w-(2w+1)\eta/4w}x^{o(\eta)},
\end{eqnarray*}
so that the corresponding contribution to (\ref{RR}) is again
$o(x^2)$.  We therefore see that (\ref{d1}) cannot hold.

We remark that this would fail for $A=B=x^{2/5}$, $C=x^{1/5}$.  It is
crucial that $C$ should not be close to $x^{1/5}$, for example, and
this is the reason for the removal of such ranges in Section \ref{S1}.

We now examine (\ref{d2}) and (\ref{d3}). Using (\ref{nx}) these become
\[R\ll_{\ep}x^{-2+2\ep}AT\{T^{1-1/w}+TC^{-1}\}B^{1+1/2w}C^2\bso
\ll_{\ep} x^{3\ep}\{x^{-1/2w}C+1\}B^{1/2w}\]
and
\begin{eqnarray*}
R^{1/4}&\ll_{\ep}& x^{-2+2\ep}AT^{3/4}\{T^{1-1/w}+TC^{-1}\}B^{1+1/2w}C^2\bso\\
&\ll_{\ep}& x^{-1/8+3\ep}\{x^{-1/2w}C+1\}B^{1/2w}
\end{eqnarray*}
  respectively.  When $w=3$ we note that
  \begin{eqnarray*}
    x^{-1/2w}CB^{1/2w}&=&x^{-1/6}CB^{1/6}\\
    &=&  x^{-1/6}(ABC)^{4/15}(C^2/A)^{11/30}(A/B)^{1/10}\\
    &\ll_{\ep}&x^{-1/6+4/15+\ep}\\
    &=& x^{1/10+\ep},
  \end{eqnarray*}
  by (\ref{ABC}).  On the other hand, if $w\ge 4$ then
\[x^{-1/2w}CB^{1/2w}\le x^{-1/2w+1/(2w-1)+1/4w+\ep}\le
x^{9/112+\ep}\ll x^{1/10},\]
by (\ref{Crng+}). Moreover
\[B^{1/2w}\le B^{1/6}\le x^{1/12+\ep}\]
for any $w\ge 3$.  We therefore see
that both (\ref{d2}) and (\ref{d3}) lead to the bound
$R\ll  x^{1/10+o(1)}$, as claimed.

\section{Factors of Length Below $x^{1/4}$}

In this section we handle the various cases in which every factor
$F_j(s)$ of $D(s)$ has length $C_j\le x^{1/4}$.  In this situation any
factor, whether of type $P_j(s)$ or $F_j(s)$, will have length at
least $z_1$ and at most
$x^{1/4-\eta}$, as shown by part (ii) of Lemma \ref{prev}.
It will be convenient to combine factors $P_j(s)$ and $F_j(s)$ of
$D(s)$ as far as possible, subject to the lengths of the resulting Dirichlet
polynomials being at
most $x^{1/4-\eta}$. Such products will no longer run over dyadic
intervals, but they will be of the form
\[Q(s)=\sum_{A<n\le 2^{3+3\nu}A}q_nn^{-s},\]
and we will refer to $A$ as being the ``length'' of $Q(s)$. Thus the
procedure described above involves multiplying any two Dirichlet
polynomials whose lengths $A_1$ and $A_2$ have
$A_1A_2\le x^{1/4-\eta}$ recursively, until no further polynomials
can be combined. We may therefore assume that $A_i\le x^{1/4-\eta}$
for any $Q_i(s)$ and that $A_iA_j>x^{1/4-\eta}$ for any two distinct
factors $Q_i(s)$ and $Q_j(s)$. If there are $m$ factors $Q_i(s)$
altogether, we deduce from (\ref{AC}) that
\beql{pbs}
2^{-3-3\nu}x\le\prod_{i=1}^m A_i\le 4x.
\eeq
We therefore see that $5\le m\le 8$. We will index the polynomials
with $A_1\ge A_2\ge\ldots$.

We begin by considering the case in which $m=5$. In view of our
ordering of the $Q_i(s)$ we will have
\[A_3A_4\le A_1A_2\le x^{1/2-2\eta}\le x^{1/2}\]
and $A_5^2\le A_3A_4$.  It follows that we can apply Proposition
\ref{pb} with 
\[A(s)=Q_3(s)Q_4(s),\;\;\; B(s)=Q_1(s)Q_2(s),\;\;\; \mbox{and}\;\;\;
C(s)=Q_5(s).\]
We then have $R\ll x^{1/10+o(1)}$ when $m=5$.

For $m=6$ we note that $A_1A_3\le x^{1/2-2\eta}\le x^{1/2}$ and
\[A_2A_4A_6\le\{A_1A_2A_3A_4A_5A_6\}^{1/2}\le(4x)^{1/2}\]
by (\ref{pbs}).
Moreover $A_5^2\le A_1A_3$ and $A_5^2\le A_2A_4A_6$. We may therefore
apply Proposition \ref{pb} with $C(s)=Q_5(s)$ and either
\[A(s)=Q_1(s)Q_3(s),\;\;\;\; B(s)=Q_2(s)Q_4(s)Q_6(s)\]
or vice-versa, depending on which of $A_1A_3$ or $A_2A_4A_6$ is
smaller.

When $m=7$ we consider two cases.  Suppose firstly that
\[A_1A_2A_6\le x^{1/2}.\]
Then $A_3A_4A_7\le A_1A_2A_6\le x^{1/2}$.
Moreover $A_5^2\le A_3A_4A_7$.  Thus we may successfully apply
Proposition \ref{pb} with 
\[A(s)=Q_3(s)Q_4(s)Q_7(s)\;\;\; B(s)=Q_1(s)Q_2(s)Q_6(s),\;\;\;
\mbox{and}\;\;\; C(s)=Q_5(s).\]
In the alternative case
we have $A_1A_2A_6\ge x^{1/2}$, whence $A_3A_4A_5A_7\le 4x^{1/2}$, by 
(\ref{pbs}). We also know that $A_1A_2\le x^{1/2-2\eta}\le x^{1/2}$.
Moreover $A_6^2\le A_1A_2$ and $A_6^2\le A_3A_4A_5A_7$. It
follows in this alternative case that we may apply Proposition
\ref{pb} with $C(s)=Q_6(s)$ and either
\[A(s)=Q_1(s)Q_2(s),\;\;\;\; B(s)=Q_3(s)Q_4(s)Q_5(s)Q_7(s)\]
or vice-versa, depending on which of $A_1A_2$ or $A_3A_4A_5A_7$ is
smaller.

There remains the case $m=8$. Here we have
\[A_2A_4A_6A_8\le\{A_1A_2A_3A_4A_5A_6A_7A_8\}^{1/2}\le(4x)^{1/2}\]
by (\ref{pbs}). Moreover $A_1\le x^{1/4-\eta}\le A_6A_8$, whence
\[A_1A_3A_5\le A_6A_8A_3A_5\le A_6A_8A_2A_4.\]
We can therefore apply Proposition \ref{pb} with
\[A(s)=Q_1(s)Q_3(s)Q_5(s),\;\;\;\; B(s)=Q_2(s)Q_4(s)Q_6(s)Q_8(s)\]
and
$C(s)=Q_7(s)$ to show that $R\ll x^{1/10+o(1)}$ in this final case.

\section{Factors of Length at Least $x^{1/4}$}

In this section we consider the case in which $D(s)$ has one or more
factors $F_j(s)$ with $C_j\ge x^{1/4}$. As in the previous section we
combine factors to produce Dirichlet polynomials $Q(s)$, but this time
we omit from the procedure any factors $F_j(s)$ for which
$C_j>x^{1/4}$. Thus any factor $Q_j(s)$ will have length
$A_j\le x^{1/4-\eta}$, and we will have $A_iA_j>x^{1/4-\eta}$ for any
distinct polynomials $Q_i(s),Q_j(s)$.

We begin by teating the case in which $D(s)$ has precisely two factors,
$F_1$ and $F_2$ say, for which $C_j>x^{1/4}$. According to part (ii) of
Lemma \ref{prev} we then have $C_j\ge x^{1/4+\eta}$.  We now write
$D(s)=F_1(s)F_2(s)H(s)$, so that the length $A$ of $H(s)$ satisfies
\[A\ll 4x/C_1C_2\ll x^{1/2-2\eta}.\]
Moreover the coefficients of $H(s)$ will have order $\bso$ by Lemma
\ref{coeffest}. 

We then deduce from (\ref{MB}) that
\[Rx \ll\bso\int_{T_1}^{2T_1}    |F_1(it)F_2(it)H(it)M(it)|dt .\]
By H\"{o}lder's inequality we therefore have
\[\{Rx\}^2\ll\bso I_1^{1/2}I_2^{1/2}
 \int_{T_1}^{2T_1}|H(it)M(it)|^2dt,\]
with
\[I_j=\int_{T_1}^{2T_1}|F_j(it)|^4dt.\]
We now apply Lemma \ref{lF4} together with Proposition \ref{P1}
to deduce that
\[{\{Rx\}^2}\ll_{\ep} \bso C_1C_2T
  \left(A^2R^2+(AT)^{\ep}\{ART+AR^{7/4}T^{3/4}\}\right),\]
for any fixed $\ep>0$. We then find that either
\[x^2\ll_{\ep} \bso C_1C_2A^2T,\]
or
\[R\ll_{\ep}C_1C_2AT^2x^{-2+2\ep},\]
or
\[R\ll_{\ep}C_1^4C_2^4A^4T^7x^{-8+8\ep}.\]

Since $C_1C_2A\le 4x$ the definition (\ref{Tdef}) of
$T$ allows us to deduce that either
\[x^{1/2}\ll_{\ep} \bs^3A,\]
or
\[R\ll_{\ep}x^{3\ep},\]
or
\[R\ll_{\ep}x^{-1/2+9\ep}.\]
The first of these is impossible by (\ref{81}), since $A\ll x^{1/2-2\eta}$,
while the other options are more than enough to give
$R\ll x^{1/10+o(1)}$.
This completes our treatment of the case in which exactly two
of the factors $F_j(s)$ have length at least $x^{1/4}$.
\bigskip

We turn now to the case in which there are three factors $F_j(s)$ with
corresponding lengths $C_j\ge x^{1/4}$, for which we use a variant of
the previous method. Writing $D(s)=F_1(s)F_2(s)F_3(s)H(s)$ we find
this time that $H(s)$ has length $A$ satisfying
\[A\le 4xC_1^{-1}C_2^{-1}C_3^{-1}\ll x^{1/4-3\eta}. \]
We then find via H\"{o}lder's inequality that
  \[Rx\ll\bso  \{I_1I_2I_3\}^{1/4}
\left\{\int_{T_1}^{2T_1}|H(it)|^4|M(it)|^4dt\right\}^{1/4},\]
  with $I_j$ as before. To estimate the remaining integral we observe that
  \[|M(it)|^4\le R^2|M(it)|^2.\]
  We may then apply Proposition
  \ref{P1} with
  \[H(it)^2=\sum_{n\le N} q_nn^{-it}.\]
  We will have $N\ll 2^{6\nu}A^2$ and $q_n\ll \bso$, by Lemma
  \ref{coeffest}. 
 A similar calculation to before then shows that either
  \[x\ll_{\ep}\bso (C_1C_2C_3)^{1/2}AT^{3/4},\]
  or
  \[R\ll_{\ep}x^{2\ep},\]
  or
  \[R\ll_{\ep}x^{-1/2+8\ep}.\]
The first of these is impossible when $A\ll x^{1/4-3\eta}$,
and the other alternatives yield $R\ll x^{1/10+o(1)}$.

Finally in this section we examine the situation in which there is
exactly one factor $F_j$ with $C_j>x^{1/4}$. Here we shall use the
following result.
\begin{lemma}\label{FASB}
  Suppose $D(s)$ factors as $F_1(s)A(s)B(s)$ with
  \[A\le 2x^{1/2-\eta},\;\;\;B\le x^{7/20}\;\;\;
  \mbox{and}\;\;\;AB\le 4x^{3/4-\eta},\]
and where $A(s)$ and $B(s)$ have no factors $F_j(s)$ of length
$C_j>x^{1/4}$.  Then $R\ll x^{1/10+o(1)}$.
\end{lemma}
From (\ref{MB}) we deduce that
\[Rx\ll\bso\int_{T_1}^{2T_1}|F_1(it)A(it)B(it)M(it)|dt,\]
  whence H\"{o}lder's inequality yields
\[Rx\ll\bso I_1^{1/4}\left\{\int_{T_1}^{2T_1}|B(it)|^4dt\right\}^{1/4}
\left\{\int_{T_1}^{2T_1}|A(it)M(it)|^2dt\right\}^{1/2}.\]
We estimate $I_1$ via Lemma \ref{lF4}, noting that
$C_1\gg 2^{-3\nu}x/AB\gg x^{1/4+\eta/2}$. To handle the second integral we
use the mean value theorem (\ref{mvdp}) coupled with
the bound $O(\bso)$ for the coefficients of $B(s)$
given  by Lemma \ref{coeffest}. The final integral can be dealt
with via Proposition \ref{P1}, again using Lemma \ref{coeffest} to
estimate the coefficients. We conclude that
\begin{eqnarray*}
  Rx&\ll& \bso\left\{C_1^2T\right\}^{1/4}
  \left\{(T+B^2)B^2\right\}^{1/4}
  \left\{\int_{T_1}^{2T_1}|A(it)M(it)|^2dt\right\}^{1/2}\\
  &\ll& \bso C_1^{1/2}(T^{1/2}B^{1/2}+T^{1/4}B)
    \{R^2A^2+x^{\ep}(RAT+R^{7/4}AT^{3/4})\}^{1/2}.
\end{eqnarray*}
Thus either
\[Rx\ll \bso C_1^{1/2}(T^{1/2}B^{1/2}+T^{1/4}B)RA
\ll \bs^2 C_1^{1/2}(x^{1/4}B^{1/2}+x^{1/8}B)RA,\]
by (\ref{Tdef}), or
\[Rx^2\ll x^{\ep+o(1)}C_1(TB+T^{1/2}B^2)AT\ll x^{2\ep}C_1(xB+x^{3/4}B^2)A,\]
or
\[Rx^8\ll x^{4\ep+o(1)}C_1^4(T^4B^4+T^2B^8)A^4T^3\ll
x^{5\ep}C_1^4(x^{7/2}B^4+x^{5/2}B^8)A^4.\]
The first alternative is impossible, since
\[C_1^{1/2}B^{1/2}A\le (4x)^{1/2}A^{1/2}\le 4x^{3/4-\eta/2}\]
and
\[C_1^{1/2}BA\le (4x)^{1/2}(AB)^{1/2}\le 4x^{7/8-\eta/2}.\]
The second option yields 
\[R\ll x^{2\ep}(1+Bx^{-1/4})\ll x^{1/10+2\ep}.\]
Finally, the third case produces
\[R\ll x^{5\ep}(x^{-1/2}+x^{-3/2}B^4)\ll 1.\]
The lemma therefore follows.
\bigskip

We are now ready to complete our treatment of the case in which
\[D(s)=F_1(s)\ldots F_k(s)P_1(s)\ldots P_t(s),\]
with $1\le k\le 3$, where $C_1\ge x^{1/4+\eta}$ and $C_j\le x^{1/4-\eta}$ for
$j\not=1$. We combine all factors other than $F_1(s)$ as far as
possible into Dirichlet polynomials $Q_i(s)$ of length
$A_i\le x^{1/4-\eta}$.  We may then write $D(s)=F_1(s)Q_1(s)\ldots Q_m(s)$
with $x^{1/4-\eta}\ge A_1\ge A_2\ge\ldots$, and 
$A_iA_j\ge x^{1/4-\eta}$ whenever $i\not=j$. We therefore
see that we must have $m\le 6$. Indeed, since $C_1\le x^{5/8}$ we must
also have $m\ge 2$. Moreover we will have
\[A_1\ldots A_m\le\frac{4x}{C_1}\le 4x^{3/4-\eta}.\]

Lemma \ref{FASB} immediately handles the cases $m=2$ and $m=3$,
by taking $B(s)=Q_1(s)$ and $A(s)=Q_2(s)$ for $m=2$, and 
$B(s)=Q_1(s)$ and $A(s)=Q_2(s)Q_3(s)$ for $m=3$. When $m=4$ the choice
$B(s)=Q_2(s)$ and $A(s)=Q_1(s)Q_3(s)Q_4(s)$ works similarly if
$A_1A_3A_4\le 2x^{1/2-\eta}$.
On the other hand, if $m=4$ and $A_1A_3A_4\ge 2x^{1/2-\eta}$ we will
have
\[C_1A_2\le 4x/(A_1A_3A_4)\le x^{1/2+O(\eta)}\]
and $A_4^2\le A_1A_3\le C_1A_2$.  Thus Proproposition \ref{pb} applies, with
$A(s)=Q_1(s)Q_3(s)$, $B(s)=F_1(s)Q_2(s)$ and $C(s)=Q_4(s)$.

When $m=5$ we apply Lemma \ref{FASB}, taking
$A(s)=Q_1(s)Q_3(s)Q_5(s)$ and $B(s)=Q_2(s)Q_4(s)$. Since
$A_1\ge A_2\ge\ldots$ we have
\[(A_1A_3A_5)^2\le C_1A_1A_2A_3A_4A_5\frac{A_1}{C_1}\le
4x\frac{x^{1/4-\eta}}{x^{1/4+\eta}}=4x^{1-2\eta},\]
whence $A\le 2x^{1/2-\eta}$.  Moreover
\[(A_2A_4)^{3/2}\le A_1A_2A_4=\frac{C_1A_1A_2A_3A_4A_5}{C_1A_3A_5}\le
\frac{4x}{x^{1/4+\eta}\cdot x^{1/4-\eta}}=4x^{1/2},\]
so that $B=A_2A_4\le x^{7/20}$ for large $x$.  The conditions of the
lemma are therefore satisfied, whence $R\ll x^{1/10+o(1)}$.

There remains the case $m=6$. Since $C_1\ge x^{1/4+\eta}$ and
$A_iA_j\ge x^{1/4-\eta}$ whenever $i\not=j$ it follows from
(\ref{AC}) that $C_1=x^{1/4+O(\eta)}$ and $A_j=x^{1/8+O(\eta)}$
for every index $j$. We may then apply Proposition \ref{pb}, with 
\[A(s)=Q_3(s)Q_4(s)Q_5(s),\;\;B(s)=F_1(s)Q_1(s)Q_2(s)\]
and $C(s)=Q_6(s)$,
again concluding that $R\ll x^{1/10+o(1)}$.
\bigskip

We have now covered all the relevant cases and have thus completed the 
proof of Theorem \ref{T2}.

\bigskip

\bigskip

Mathematical Institute,

Radcliffe Observatory Quarter,

Woodstock Road,

Oxford

OX2 6GG

UK

\bigskip

{\tt rhb@maths.ox.ac.uk}


\begin{thebibliography}{99}

\bibitem{BCR}S. Bettin, V. Chandee, and M. Radziwi\l\l, The mean
  square of the product of the Riemann zeta-function with Dirichlet
  polynomials, {\em J. Reine angew. Math.}, 729 (2017), 51--79.

\bibitem{HBI}D.R. Heath-Brown, The differences between consecutive
  primes, {\em J. London Math. Soc.}, 18 (1978), 7--13.
  
\bibitem{HBIII}  D.R. Heath-Brown, The differences between
consecutive primes, III, {\em J. London Math. Soc.} 20, (1979),
177--178.

\bibitem{hbi}D.R. Heath-Brown, Prime numbers in short intervals and a
  generalized Vaughan identity, {\em Can. J. Math.}, 34 (1982),
  1365--1377.
\bibitem{hbg}D.R. Heath-Brown,   The number of primes in a short interval,
{\em J. Reine angew. Math.}, 389 (1988), 22--63.
  
\bibitem{RHBIV}  D.R. Heath-Brown, The differences between consecutive
  primes, IV,  {\em A tribute to Paul Erd\H{o}s},
  (Cambridge University Press, 1990), 277--287.

\bibitem{smooth}D.R. Heath-Brown, The differences between consecutive
  smooth numbers, {\em Acta Arithmetica,}  184 (2018), 267-285.

\bibitem{hux}M.N. Huxley, On the difference between consecutive
  primes, {\em Invent. Math.}, 15 (1972), 164--170.

\bibitem{1114}M. Jutila, Zero density estimates for $L$-functions,
  {\em Acta Arith.}, 32 (1977), 55--62.

  
\bibitem{Mat}K. Matom\"{a}ki, Large differences between consecutive 
  primes. {\em Q. J. Math.} 58 (2007), 489--518.

\bibitem{mont}H.L. Montgomery, {\em Topics in multiplicative number
  theory}, Lecture Notes in Math. 227, (Springer, Berlin-Heidelberg-
  New York), 1971.


\bibitem{Peck}A.S. Peck, Differences between consecutive primes,
  {\em Proc. London Math. Soc.}, 76 (1998), 33--69.


\bibitem{Selb}A. Selberg, On the normal density of primes in small
intervals and the difference between consecutive primes, {\em Arch. Math. 
  Naturvid.}, t. 47, (1943), no. 6, 87--105.

\bibitem{shiu}P. Shiu, A Brun--Titchmarsh theorem for multiplicative
  functions, {\em J. Reine angew. Math.}, 313 (1980), 161--170.

\bibitem{titch}E.C. Titchmarsh, {\em The theory of the Riemann 
zeta-function}. Second edition. (Oxford University Press, New York, 1986).

\bibitem{W}D. Wolke, Grosse Differenzen aufeinanderfolgender
Primzahlen, {\em Math. Ann.}, 218 (1975), 269--271.

\bibitem{YG}G. Yu, The differences between consecutive primes,
  {\em Bull. London Math. Soc}. 28 (1996), 242--248.





\end{thebibliography}
\end{document}